\documentclass{article}
\usepackage[utf8]{inputenc}

\usepackage[utf8]{inputenc}
\usepackage{amsthm}
\usepackage{amsmath}
\usepackage{amssymb}
\usepackage{verbatim}
\usepackage{float} 
\usepackage[all]{xy}
\usepackage{graphicx}
\usepackage[bottom=3cm,top=3cm,left=3cm,right=3cm]{geometry}
\newtheorem{teo}{Theorem}
\newtheorem{supo}{Assumption}
\newtheorem{prop}{Proposition}

\begin{document}

\newpage
\begin{center}

\LARGE{\textbf{Misspecification Tests on Models of Random Graphs}}

\Large{Denise Duarte*  and
Rafael Honório Pereira Alves*}
\end{center}

*Departamento de Estatística da Universidade Federal de Minas Gerais

\thispagestyle{empty}
\pagenumbering{gobble}

\begin{center}
  \large{\textbf{Abstract}}
\end{center}

In recent years, there has been a great interest in random graph models to model complex networks in the most diverse areas such as Social Sciences, Physics, Biology, Economics, Ecology and Computer Science.

 A class of models that have been widely used are the exponential random graph (ERG) models, which form a comprehensive family of models that include independent and dyadic edge models, Markov random graphs, and many other graph distributions, in addition to allow the inclusion of covariates that can lead to a better fit of the model.

Another increasingly popular class of models in statistical network analysis are stochastic block models (SBMs). They can be used for the purpose of grouping nodes into communities or discovering and analyzing a latent structure of a network. The stochastic block model is a generative model for random graphs that tends to produce graphs containing subsets of nodes characterized by being connected to each other, called communities.

 Many researchers from various areas have been using computational tools to adjust these models without, however, analyzing their suitability for the data of the networks they are studying. The complexity involved in the estimation process and in the goodness-of-fit verification methodologies for these models can be factors that make the analysis of adequacy difficult and a possible discard of one model in favor of another.
 And it is clear that the results obtained through an inappropriate model can lead the researcher to very wrong conclusions about the phenomenon studied.

The purpose of this work is to present a simple methodology, based on Hypothesis Tests, to verify if there is a model specification error for these two cases widely used in the literature to represent complex networks: the ERGM and the SBM. We believe that this tool can be very useful for those who want to use these models in a more careful way, verifying beforehand if the models are suitable for the data under study.

\vspace{.5cm}


\pagenumbering{arabic}

\section{Introduction}
The understanding of interaction mechanisms in complex real-world networks through random graph models is a topic that has gained a lot of attention in the literature in recent years Newman(2001).

Practical applications of random graphs are found in all areas where complex networks need to be modeled, some examples are Social Science, Physics, Biology, Economics, Ecology and Computer Science. Economic or social interactions are also often organized into complex network structures. Similar phenomena are observed in communication networks such as the internet or in traffic flow. In current problems in Biosciences, protein networks in the cell are important examples, as well as molecular networks in the genome. On larger scales, networks of cells are found as in neural networks, up to the scale of organisms in ecological food webs. Many random graph models try to mirror the different types of complex networks found in different areas. For a more comprehensive description of the theme and applications, we recommend reading Newman (2010) and Reuven and Shlomo (2010).

The random graph of Erd\"{o}s and R\'{e}nyi( Erd\H{o}s, 1959) is one of the best studied network models, however, for real world networks such as social networks, the Internet or biological networks, it is not a good model as some basic properties do not fit well. Complex networks tend to have non-trivial topological characteristics that differ from random graphs of Erd\"{o}s  Renyi, such as heavy tail in the degree distribution, high clustering coefficient, hierarchical structures and average length of short paths. Two models of random graphs widely used to model natural phenomena that try to capture these characteristics are the scale-free networks, proposed by Barab\'{a}si and Albert (1998) and the models of Small-World introduced by Watts and Strogats (1989).

In recent years, there has been a great interest in exponential random graph models to model networks, especially social networks. (Frank and Strauss, 1986, Frank, 1991, Wasserman and Pat- tison,  1996; see also  Pattison  and  Wasserman,  1999,  Robins  et  al.,  1999).   The exponential random graph model class is a comprehensive family of models that includes independent and dyadic edge models, the random graphs of Markov  de  Frank  and  Strauss  (1986and many other graph distributions, in addition to allowing the inclusion of covariates that can lead to a better fit of the model. The estimation of the parameters of this model is done through computational methods and is implemented in several programming languages, in software R, for example, has the package called ERGM implemented (Hunter et al, 2009).

Another increasingly popular class of models in the statistical analysis of networks are the stochastic block models (SBMs) introduced by Holland  et  al  (1983).   They can be used for the purpose of grouping nodes into communities or discovering and analyzing a latent structure of a network. There has been rapid development in the theme of clustering based on random graph models in the last ten years. Citing just a few, we have the two works by Abbe and Sandom ( 2015)  and that of Karrer  e  Newman  (2011).  The stochastic block model is a generative model for random graphs that tends to produce graphs containing subsets of nodes characterized by being connected to each other, called communities. This model is hierarchical in the sense that we first sort the nodes that belong to each group, then we sort the edges between the nodes. Each of these draws is governed by a law that depends on specific parameters. Due to the complexity of the likelihood of this two-stage process, the estimation of parameters of a SBM is also based on computational methods that are implemented in several languages.   For example, in the R software, the Stochastic Block Model (SBM) package implemented by Leger (2016) is implemented. A detailed review of SBMs can be found in Lee and Wilkinson (2019) and a very interesting extension applied to large networks is presented in Peixoto (2014).

Many researchers from various areas have been using computational tools to fit models without, however, making a critical analysis regarding its suitability for the data being analyzed, especially in cases where the chosen model is an ERG or a SBM. For the ERG we cite, for example, the verification procedures based on resampling presented in  Kolaczyk  e  Csárdi  (2014).   In the case of SBM, one of the model-fitting verification methodologies are based on the Akaike Information Criterion, AIC (Akaike, 1973 and 1974)) and are implemented in the SBM package of R. The complexity involved in the estimation process and in the quality-of-fit verification methodologies for these models may be factors that hinder the analysis of adequacy and a possible discard of one model in favor of another. It is clear that the results obtained through an unsuitable model can lead the researcher to very wrong conclusions about the phenomenon studied.

The purpose of this work is to present a simple methodology, based on Hypothesis Tests, to verify if there is a model specification error for these two cases widely used in the literature to represent complex networks: the ERG and the SBM. We believe that this tool can be very useful for those who want to use these models in a more careful way, verifying beforehand if the models are suitable for the data under study. We highlight the advantage of this type of methodology in relation to the selection processes of models of the AIC or BIC type (Schwarz,1978) and others with the same proposal, because in these tools we do not have an indication as to the suitability of the models, but a selection from among the proposed candidates, which may all be unsuitable.

We will derive the tests taking into account the likelihood of these two models, using maximum likelihood estimators already proposed in the literature, for the ERG we quote Schmid and Desmarais (2017) and for the SBM we quote Celisse, Daudin and Pierre (2012). We will closely follow the work presented by White (1982), where a model fit verification test is proposed for random variables whose densities satisfy a general set of regularity conditions.We will show that the ERG and SBM models satisfy these conditions and we will explain the statistics that must be calculated to perform the tests on each of these models. We will develop a package in R with the implementation of tests for the two models and we will present a simulation study showing that the proposed tests achieve the desired purposes.

\section{Main Definitions}
\subsection{Some models for random graphs}
\subsubsection{Erdos Rényi Random Graphs}

A random graph $G$ is a random variable that takes values in a family of graphs $\cal{G}$. The study of such graphs dates back to the 1950s, when Paul Erdos and Alfréd Rényi derived a series of results on random graphs. The graphs we refer to here are all labeled, that is, the vertices are distinct (for example, there are $\binom{n}{2}$ graphs with $n$ vertices and exactly one edge). 

There are two ways to define the Erd\H{o}s-Rényi random graph model: the model denoted by $G(n,M)$, represents a graph chosen at random from the collection of all graphs with $n$ nodes and $M$ edges. For example, in the model $G(3,2)$ each of the three possibilities of graphs with three vertices and two edges are included with a probability of $\frac{1}{3}$.

And the Erdos-Rényi model denoted by $G(n, p)$, represents a graph with $n$ vertices in which: for each pair $u, v$ of vertices of $G = G(n, p)$ , the edge $uv$ is random (exists or does not exist), regardless of any other edge. A Bernoulli independent random variable with parameter $p$ is used to decide on the presence of the edge $uv$; the edge $uv$ is in the graph if and only if the variable results in success. By the independence between the edges, the probability of a graph $G$, with $n$ vertices and $M$ edges, is 
\begin{equation}
P(G) = p^M(1 - p)^{\binom{n}{2}-M}. \label{er}
\end{equation}

Several properties of Erdos-Rényi graphs are now well known. This model is sometimes used as a reference for evaluating others, representing the role of a network with uniform random behavior.(See Frank e Strauss (1986))

Several other random graph models are studied in the literature. They attend the need for more complicated and realistic models(than the ER model) that describe the real networks observed in practice. Some characteristics considered typical of real networks are the presence of a large number of vertices $n$, few edges (e(G) = O(n)), small diameter ($diam(G) = O(\log n) $, ie, two vertices taken at random are connected by a short path), degrees of vertices distributed according to a power law (the number of vertices
with degree $k$ is proportional to $k^{-\beta}$, for some constant $\beta$), and clustering effect (clustering, or transitivity of connections: vertices with common neighborhoods are more likely to be connected) . The number of publications and studies of complex networks and random graphs is large. For an introduction to the area, we cite Chatterjee and Diaconis (2011) and Robins, Pattison and Kalish (2007).

\subsubsection{Exponential Random Graphs (ERG) }

The ERG model, also known as the $(p*)$model, was first proposed by Holland and Leinhardt (1983), and was built on the statistical foundations established by Besag (1974). These models constitute a family of statistical models that have been widely used to model social networks. The importance of this model lies in its ability to represent the social structural effects commonly observed in many human social networks, including general effects based on the degree of each vertex, such as reciprocity and transitivity, or even, activity based attributes and popularity effects.

Substantial developments were made by Frank and Strauss (1986), and continued to be made by other authors throughout the 1990s. And they are still the object of study by several researchers, always in an attempt to make the ERG an even better model, more applicable to real network data. A detailed review of the subject is presented in Wasserman and Pattison (1996). The model proposed by Besag (1979), in the context of Spatial Statistics, in turn, is centered on the Hammerley-Clifford Theorem (1971), born in Statistical Physics, which shows that a probabilistic model for graphs, with a certain dependence structure must necessarily belong to an exponential family. More than that, This important Theorem also shows which network information should be used to calculate the probabilities of a given configuration.

The problem that arises when considering dependency structures between vertices, even though they are considerably simple structures like those considered in the article mentioned above, is that the number of parameters to be estimated in the model is very large, as we will detail later in this text. The inferences proposed for the parameters of this model are based on pseudo maximum likelihood, using an analogy with the logistic regression model, which makes the inferences unreliable. In an attempt to overcome these and other problems, in the 2000s new specifications on the vertex dependence structure were proposed by researchers, such as those found in the series of articles published by Snijders (1997). The main idea in these works is to reduce the dimension of the vector of parameters of the models and preserve the characteristics of the network. ERG models also appear that allow exogenous characteristics of the network to be used to model the probability of occurrence of a configuration, thus improving inferences for the network. An alternative ERG model to those proposed by the articles by Snijders and collaborators, but with a similar objective in relation to the decrease in the number of parameters, was proposed by Hunter and Handcock (2009).

A random graph model is an ERG if, for every graph $G \in \cal{G}$, we can express the probability of its occurrence by
$$P(G; \theta)=\exp\left( \sum_{i=1}^{n}\theta_iT_i(G)-\varphi(\theta)\right)=\exp\left( \sum_{i=1}^{n}\theta_iT_i(G)\right)/z(\theta).$$
where $\theta=(\theta_1, \ldots , \theta_n)$ is a vector of known real parameters and $T_i(G)$ are functions of $G$ (such as the number of edges, triangles, stars, circuits, etc. .); that is, if the distribution over the graph space is a member of the exponential family of distributions. As in the Erdos-Rényi model, we consider the number $n$ of vertices of $G$ fixed. The factor $e^{-\varphi(\theta)}=z^{-1}\left(\theta\right)$ is sometimes called the normalization constant.

Two difficulties in using exponential random graphs are the estimation of $z(\theta)$ and the fact that very different values of $\theta$ give rise to essentially equal distributions in the graph space (See Chatterjee and Diaconis, 2011) .
\begin{itemize}
\item{ Independent edge case}\\
Assume that the connections between the vertices occur independently of each other, that is, that there is no dependence within the network. In this case, the function $T_i(G)$ becomes only the pointer of the edge $Y_{ij}$ of the adjacency matrix $Y$ of $G$, in this way, the general ERG model is greatly simplified, since the parameters of the model reduce to the binding coefficients $ij$ and the ERG reduces to the Erdös-Rènyi model.

\item{Markov Graph Model}\\
Following the work of Besag (1974) in the area of Spatial Statistics, Frank and Strauss (1986) proposed a Markov dependence on a Graph, postulating that a possible link from i to j is assumed to be dependent on any other possible link involving i or j, even though all other connections in the network are fixed. Markov dependence implies that two possible edges of a network are conditionally independent unless they share a common vertex. They showed that this assumption resulted in models for undirected graphs that involve parameters associated with simple network statistics such as number of edges, star-shaped and triangle-shaped structures. In this model, two vertices are considered neighbors if they share an edge. A subset of the set of vertices, V, where all elements are neighbors is called a clique.

We note that all model specifications involve statistics that are just functions of the y network itself, only endogenous effects are considered. Even so, it is natural to expect that the probability of a connection between two vertices may also depend on characteristics, attributes of the vertices themselves. So, allowing the incorporation of exogenous effects can lead to more accurate inferences about the network. We can incorporate attributes that were measured at the verties, in the form of additional statistics in the function within the exponential.

\end{itemize}

\subsubsection{Stochastic Block Models (SBM) }

When analyzing complex networks, a basic task in the area of community detection (or clustering) is to partition the vertices of a graph into clusters that are more densely connected. More generally, community structures can also refer to groups of vertices that connect similarly to the rest of the graph, without necessarily having a higher internal density. In the most general context, community detection refers to the problem of inferring similarity relationships between items in a network by observing their local interactions.

Community detection is one of the central problems in networking and data science. The Stochastic Block Model(SBM) has been widely used as a canonical model to study these issues.

Let's define the Stochastic Block Model: A graph $G$ in the set of vertices $v(G)$ can be represented by its adjacency matrix $Y=\left\{ Y_{ij} \right\}_{1 \leq i \neq j \leq n}$, where
$$Y_{ij}=\left\{ \begin{array}{cc}
1 & \text{ if there is an edge between the vertices $i$ and $j$} \\
0 & \text{ otherwise,}
\end{array}\right.$$
where $Y_{ii}=0$ for all $i$, that is, there is no connection of the vertex with itself.

We consider a graph whose vertices belong to different $m$ categories. These categories we will call blocks. Let $X=\left( X_i \right)^{n}_{i=1}$, where $X_i=k$, if vertex $i$ belongs to block $k$, for all $i \in \left\{1, \ldots,n \right\}$ and $k \in \left\{1, \ldots,m \right\}$. Then the block graph can be represented by $(Y,X)$, $X$ is called the block structure of the graph $G$.

For a random block graph, the number of vertices $n$ is fixed, but the adjacency matrix $Y$ and the block structure $X$ are random.

Let the vertex set be $\left\{1 , \ldots, n \right\}$ and the following conditions:
\begin{itemize}
\item[$1.$]$Y_{ij}=Y_{ji}$ and $Y_{ii}=0$.
\item[$2.$]There is a partition of the $n$ vertices into $m$ blocks such that for all $i,j,h$ with $i \neq j \neq h$, if $i$ and $h$ belong to the same block, then $Y_{ij}$ and $Y_{hj}$ are identically distributed.
\end{itemize}
The stochastic block model is a probability distribution family of a block graph $G$ with vertex set $\left\{1 , \ldots, n \right\}$ and block set $\left\{1 , \ldots, m \right\}$, defined as follows:
\begin{itemize}
\item[$1.$]Parameters are the vector $\theta=(\theta_1, \ldots, \theta_m)$, of the block probabilities and the matrix $\eta=\left( \eta_{kl}\right) _{1 \leq k \leq l \leq m}$, from the probabilities of the block-dependent edges.
\item[$2.$]The vector of blocks consists of the $\left(X_i\right)^{n}_{i=1}$ independent and identically distributed random variables, where $P(X_i=k)= \theta_k$, for $k=1 \ \ldots, m$.
\item[$3.$]Conditional to the block of vertex $X_i$, edges $Y_{ij}$ are independent with $Y_{ij} \sim \text{Bernoulli}(\displaystyle \eta_{X_i, X_j}) $.
\end{itemize}
If $(X,Y)$ represents a block graph $G$, the probability function is given by:
$$P( \theta, \eta ; X,Y)=\theta_{1}^{n_1} \cdots \theta_{m}^{n_m}\prod_{1 \leq k \leq l \leq m}\eta_{kl}^{e_{kl}}(1-\eta_{kl})^{n_{kl}-e_{kl}},$$
where $n_k=\sum_{i=1}^{n}I(X_i=k)$ denotes the number of vertices of $G$ that belong to block $k$,
$$e_{kl}=\sum_{1 \leq i \neq j \leq n}Y_{ij}I(x_i=k)I(x_j=l).$$
denotes the number of edges of $G$ that have a vertex in block $k$ and a vertex in block $j$, and
$$n_{kl}=\left\{ \begin{array}{cc}
n_k n_l & \text{ if } k \neq l \\
\binom{n_k}{2} & \text{ if } k=l,
\end{array}\right.$$

We will also denote $s=\sum_{1 \leq i \leq j \leq n} Y_{ij}$ the total number of edges.

The conditional distribution of the block graph given the block vector $\left( X_i \right)^{n}_{i=1}$ is a stochastic block model with independent edges where the blocks are a function of the parameters. In general, the number of model parameters tends to infinity together with $n$, which makes their estimation difficult. Several stochastic properties of stochastic block models are studied in the literature. (See Snijders (1997) and Celisse, Daudin and Pierre (2012))

\subsection{Maximum Likelihood Estimation}

\subsubsection{Maximum Likelihood Estimation on ERG}
Suppose $X_1, \ldots, X_n$ are independent and identically distributed random variables
following distribution $f(\cdot| \theta)$. Given the observed values $x_1, \ldots, x_n$ we can build the likelihood function:
$$\mathcal{L}(\theta; x_1, \ldots, x_n)=\prod_{i=1}^{n}f(x_i|\theta).$$
This function is the joint density of $x_1, \ldots, x_n$, but as a function of $\theta$.
Let $\hat{\theta}$ be the value of $\theta$ corresponding to the global maximum of the function, $\hat{\theta}$ is called the model's maximum likelihood estimator.
An easier way to find $\hat{\theta}$ is to use the log-likelihood function.
This function has the same maximum likelihood estimator $\hat{\theta}$ as likelihood function. The log-likelihood function is defined by:
$$\textit{l}(\theta; x_1, \ldots, x_n)=\log \left\{ \prod_{i=1}^{n}f(x_i|\theta) \right\}=\sum_ {i=1}^n \log \left[f(x_i| \theta)\right].$$
In the case of the ERG we have that
$$\mathcal{L}(\theta; T)=\prod_{i=1}^{n} \frac{\exp(\theta_i T_i(G))}{z(\theta)},$$
implies in
$$\textit{l}(\theta; T_1, \ldots, T_n)=\sum_{i=1}^n \theta_i T_i(G)- \log \left\{ \sum_{y \in \Omega} \exp\left[ \theta_i T_i(y)\right]\right\}$$
where $\Omega$ is the set of all possible graphs of $n$ vertices. 

In the simplest cases, for example in the one-parameter model, direct maximization is easy to obtain:
  $$\frac{e^{\hat{\theta}}}{1+e^{\hat{\theta}}}=\frac{\sum_{t=1}^{\binom{n}{2 }}U_t}{\binom{n}{2}},$$
 
  where $U_t$ are the $\binom{n}{2}$ random variables that represent the existence of the edges of the graph and will be discussed in detail throughout this thesis.
 
  When we are dealing with very large networks, it is very difficult to differentiate the second term of the equation and the computational complexity increases with the number of enough statistics used. In both cases the estimator obtained is consistent.

Some methods are used in the literature, for example the pseudo-likelihood estimation method or Markov Chain's Monte Carlo method (MCMC) which is most used and implemented today. (See Corander and Dahmstrom (1998))

\subsubsection{Maximum Likelihood Estimation on SBM}
In the case of SBM, we have that the likelihood function is given by
$$\mathcal{L}(\theta, \eta; x, y)=\prod_{i=1}^{n}\prod_{j=1}^{n} \left\{\prod_{l= 1}^{m}\theta_1^{\left(\textit{I}_{x_i=l}\right)}\prod_{1 \leq k \leq l \leq m}\eta_{kl}^{e_ {kl}}(1-\eta_{kl})^{n_{kl}-e_{kl}}\right\},$$
implies in
$$\textit{l}(\theta, \eta; x, y)=\sum_{i=1}^{n}\sum_{j=1}^{n}\log \left( \left\{ e^{\sum_{k=1}^{l}\sum_{l=1}^{m}\left\{ e_{kl} \log( \eta{kl})+(1-e_{kl} )\log(1-\eta_{kl})\right\}}\prod_{l=1}^{m}\theta_1^{\left(\textit{I}_{x_i=l}\right)} \right\} \right).$$

This function is not easy to maximize, firstly because $Y_{ij}$ is only conditionally independent of $X_i$ and $X_j$, secondly because the number of random variables in the expression \begin{equation} \sum_{ k=1}^{l}\sum_{l=1}^{m}\left\{ e_{kl} \log( \eta{kl})+(1-e_{kl})\log(1- \eta_{kl})\right\}\label{lv}\end{equation} is $\binom{n}{2}$, much greater than the number of vertices $n$.

Several techniques are used in this case. A direct maximization can be done computationally by applying a log-likelihood transformation that passes the expression \ref{lv} from a problem with $\binom{n}{2}$ variables to a problem with $n$ variables. The Expectation Maximization(EM) algorithm can be applied, for that, it is necessary to take the vector $X$ as the missing data vector and maximize the log-likelihood expectation restricted to $P(Y,X, \theta,\eta |Y, \theta, \eta)$. One can still use MCMC methods. All methods generate consistent estimators, however in practice the estimation can only be done for graphs with reasonably small $n$ and $m$ and generate generally unstable estimates of $\eta$. (See Snidjers (1997))

An estimation method using EM and variational techniques, called Variational EM, was proposed and has been the most used in the study of SBM. The method can still generate unstable estimates for $\eta$, but it is quite useful in practice, as it can estimate the parameters for graphs with large $n$ and/or $m$ and the estimators obtained are consistent. (See Celisse, Daudin and Pierre(2012)). The package \textit{blockmodels} of the R software uses the Variational EM to estimate the parameters of the SBM (Leger, 2016).

\subsection{Quasi-Likelihood Estimation for Misspecification}
Since Fisher postulated the maximum likelihood method in the 1920s, the method has become one of the most important tools for estimation and inference available to statisticians.

A fundamental assumption underlying the classical results on the properties of the maximum likelihood estimator is that the stochastic law that determines the behavior of the phenomena investigated (the true structure) is known within a specified parametric family of probability distributions (the model). In other words, the probability model is considered to be correctly specified. In many (if not most) circumstances, one may not have complete confidence that this is so. If the probability model is not assumed to be correctly specified, it is natural to ask what happens to the properties of the maximum likelihood estimator. Does it still converge to some limit asymptotically, and does that limit have any meaning? If the estimator is somehow consistent, is it also asymptotically normal? Does the estimator have properties that can be used to decide whether or not the specified family of probability distributions contains the true structure? We will provide the answers to these questions.

Under some conditions, the maximum quasi-likelihood estimator (QMLE) is a natural estimator for the parameters that minimize the Kullback-Leibler information criterion, so the maximum likelihood estimator converges to a well-defined limit, even when the probability is not specified correctly. An interesting feature of this result is that, with the wrong specification, the asymptotic covariance matrix of the QMLE is no longer equal to the inverse of the Fisher information matrix. However, the covariance matrix can be estimated consistently and, as expected, simplifies the familiar form in the absence of specification errors. This property is exploited to produce a new test for specification errors, applicable to a wide range of problems, in this work we will extend to ERG and SBM. (See White (1982))

The quasi-log-likelihood function of the sample is defined as follows:
$$L_n(U,\theta)\equiv n^{-1}\sum_{t=1}^{n}\log f(U_t, \theta)$$
and we define the quasi-maximum likelihood estimator (QMLE) as the vector of parameters $\hat{\theta_n}$ which is the solution of the equation:
$$\hat{\theta_n}=\max_{ \theta \in \Theta}L_n(U,\theta)$$

Next, we present some assumptions that will be necessary to show the existence and convergence results of QMLE presented in White (1982). The proofs of the Theorems stated in this section can also be found in this article.

\begin{supo} The independent random vectors $1 \times M$, $U_t$ with $t= 1,\ldots, n$, have a common joint distribution function $H \in \Omega$ , a measurable Euclidean space, with measurable Radon-Nikodym density $h = dH / d \nu$. \label{supo1}
\end{supo}

As $H$ is unknown a priori, we choose a family of distribution functions that may or may not contain the true structure, $H$. It is usually easy to choose this family to satisfy the next assumption.

\begin{supo}The family of distribution functions $F (u, \theta)$ has Radon-Nikodym densities $f (u, \theta) = dF(u, \theta) / d\nu$ that are measurable in $u$ for each $\theta \in \Theta$, a compact subset of a $p$-dimensional Euclidean space and continuous in $\theta $ for each $u \in \Omega$. \label{supo2}
\end{supo}

\begin{teo}
Given the assumptions \ref{supo1} and \ref{supo2}, for all $n$ there is a measurable QMLE $\hat{\theta_n}$.\label{teo1}
\end{teo}

Once the existence of a QMLE is assured, we move on to examining its properties. When $F$ contains the true structure $H$ (that is, $H (u)=F (u, \theta_0)$ for some $\theta_0 \in \Theta) $ the general theory of maximum likelihood estimators guarantees that the MLE is consistent for $\theta_0$ under proper regularity conditions. However, without this restriction, he observed that since $L_n(U,\theta)$ is a natural estimator of $E (\log f (U_t,\theta))$, $\hat{\theta_n}$ is a natural estimator of $\theta_{*}$, the parameter vector that minimizes the Kullback-Leibler Information Criterion (KLIC),
$$I(h:f, \theta)\equiv E\left(\log\left[\frac{ h(U_t)}{f(U_t, \theta)}\right]\right).$$
Here expectations are taken with respect to the true distribution. Therefore,
$$I(h:f, \theta)\equiv \int \log(h(u))dH(u)- \int \log(f(u,\theta))dH(u).$$
The opposite of $I(h:f, \theta)$ is called the entropy of the distribution $H (u)$ with respect to $F (u, \theta)$. Intuitively, $I(h:f, \theta)$ measures our ignorance of the true structure. For $\hat{\theta_n}$ to be a natural estimator of $\theta_{*}$, we impose the following condition
\begin{supo}
\label{supo3}
\begin{itemize}
\item[a)]$E(\log( h(U_t))$ exists and $|\log f(U_t, \theta)|\leq m(u)$ for all $\theta \in \Theta$, where m is integrable with respect to $H$
\item[b)]$I(h:f, \theta)$ has a unique minimum at $\theta_{*} \in \Theta$.
\end{itemize}
\end{supo}

The assumption \ref{supo3} ensures that the KLIC is well-defined.

\begin{teo}
Given the assumptions \ref{supo1} to \ref{supo3}, $\hat{\theta_n} \rightarrow \theta_{*}$
when $n\rightarrow \infty$ for almost every sequence $(U_t)$. \\($\hat{\theta_n} \overset{q.c.}{\rightarrow} \theta_{*}$)\label{teo2}.
\end{teo}

In other words, the QMLE is generally a consistent estimator for vector parameters.

The next step is to show an asymptotic normality of the QMLE and, for that, it is necessary to define some auxiliary matrices, when partial derivatives exist:
$$A_n(\theta)=n^{-1}\sum_{t=1}^{n}\left[\frac{\partial^2 \log f(U_t, \theta)}{\partial \theta_i \partial \theta_j}\right],$$
$$B_n(\theta)=n^{-1}\sum_{t=1}^{n}\left[\frac{\partial \log f(U_t, \theta)}{\partial \theta_i}. \frac{\partial \log f(U_t, \theta)}{\partial \theta_j}\right].$$
And we consider the expectations,
$$A(\theta)=E\left(\frac{\partial^2 \log f(U_t, \theta)}{\partial \theta_i\partial \theta_j}\right),$$
$$B(\theta)=E\left(\frac{\partial \log f(U_t, \theta)}{\partial \theta_i} .\frac{ \partial \log f(U_t, \theta)}{ \partial \theta_j}\right).$$
When the proper inverses exist, we define:
$$C_n(\theta)=A_n(\theta)^{-1}B_n(\theta)A_n(\theta)^{-1},$$
$$C(\theta)=A(\theta)^{-1}B(\theta)A(\theta)^{-1}.$$

\begin{supo}$\partial \log f(u, \theta)/\partial \theta_i$, $i=1, \ldots , p$, are measurable functions of $u$ for each $\theta$ in $\Theta$ and continuous functions of $\theta$ for each $u$ in $\Omega$.
\label{supo4}
\end{supo}

\begin{supo}$|\partial^2 \log f(u, \theta) / \partial \theta_i \partial \theta_j|$ e $|\partial \log f(u, \theta) / \partial \theta_i . \partial \log f(u, \theta) /\partial \theta_j|$, $i,j=1 , \ldots, p$ are dominated by integrable functions with respect to $H$ for all $u$ in $ \Omega$ and $\theta$ in $\Theta$.
\label{supo5}
\end{supo}

\begin{supo}
\begin{itemize}
\item[a)]$\theta_{*}$ is an interior point of $\Theta$,
\item[b)]$B(\theta_{*})$ is non-singular,
\item[c)]$\theta_{*}$ is a regular point of $A(\theta).$
\end{itemize}
\label{supo6}
\end{supo}

The assumption \ref{supo4} ensures that the first two derivatives with respect to $\theta$ exist, that these derivatives are measurable with respect to $\theta$ follows from the assumption \ref{supo2}, since the derivative can be considered as the limit of a sequence of measurables. These conditions allow us to apply a mean value theorem to random functions. The premise \ref{supo5} guarantees that the derivatives are dominated by integrable functions with respect to $H$, which guarantees that $A(\theta)$ and $B(\theta)$ are continuous on $\theta$ and that we can apply a law of large numbers to $A_n(\theta)$ and $B_n(\theta)$ . On the assumption \ref{supo6}, we define a regular point of the matrix $A(\theta)$ as a value for $\theta$ such that $A(\theta)$ has \textit{rank} constant in some open neighborhood of $\theta$.

With these additional assumptions, the following Theorem guarantees that the QMLE has an asymptotically normal distribution.
\begin{teo}\textbf{(Asymptotic Normality)}Given the assumptions from \ref{supo1} to \ref{supo6}:
$$\sqrt{n}\left(\hat{\theta_n}-\theta_{*}\right) \overset{A}{\sim} N(0,C(\theta_{*})).$$
Also $C(\hat{\theta_n})\overset{q.c.}{\sim}C(\theta_{*})$, element by element.
\label{teo3}
\end{teo}

We have asymptotic normality since
\begin{equation}\int \frac{\partial^2 \log f(u, \theta)}{\partial \theta_i\partial \theta_j} . f(u,\theta)d\nu=-\int \frac{\partial \log f(u, \theta)}{\partial \theta_i} .\frac{ \partial \log f(u, \theta) }{\partial \theta_j} . f(u,\theta)d\nu
\label{hi}
\end{equation}
The equation \ref{hi} is the familiar equality in maximum likelihood theory that ensures the equivalence of the Hessian(left side) and inverse Fisher Information Matrix(right side). In the present case, this equivalence will generally not be valid. However, when the model is specified correctly and the next assumption is valid, we get an information matrix equivalence result.

The following assumptions are necessary to state the Theorem \ref{teo4} about the Fisher Information associated with $\theta$.

\begin{supo}$\partial[\partial \log f(u, \theta)/\partial \theta_i. f(u,\theta)]/\partial\theta_j$, $i,j=1, \ldots, p$, are dominated by integrable functions with respect to $\nu$ for all $\theta in \Theta$, and the minimal support of $f(u,\theta)$ does not depend on $\theta$.
\label{supo7}
\end{supo}

Together, the given conditions from \ref{supo1} to \ref{supo7} and $h(u)=f(u,\theta_0)$ for some $\theta_0$ in $\Theta$, can be considered as the \textit{usual maximum likelihood regularity conditions}, as they ensure that all familiar results hold.

\begin{teo}\textbf{(Information Matrix)}Given the assumptions from \ref{supo1} to \ref{supo7}, if $g(u)=f(u,\theta_0)$ for some $\theta_0 $ into $\Theta$, then $\theta_{*}=\theta_0$ and $A(\theta_0)=-B(\theta_0)$ , hence $C(\theta_0)=-A(\theta_0)^{ -1}=B(\theta_0)^{-1}$, where $-A(\theta_0)$ is Fisher's information matrix.
\label{teo4}
\end{teo}

The Theorem \ref{teo4} essentially says that when the model is specified correctly, the information matrix can be expressed in Hessian form, $-A(\theta_0)$ or in product form, $B(\theta_0) $. Equivalently, $A(\theta_0)+B(\theta_0)=0$. When this equality fails, it follows that the model is misspecified, and this misspecification can have serious consequences when standard inferential techniques are applied. So $A(\theta_{*})+B(\theta_{*})$ is a useful indicator for misspecification!

The matrix $A(\theta_{*})+B(\theta_{*})$ is not observable, but can be consistently estimated by $A_n(\hat{\theta_{n}})+B(\hat{ \theta_{n}})$. To obtain a test statistic, we consider the asymptotic distribution of the elements of $\sqrt{n}(A_n(\hat{\theta_{n}})+B(\hat{\theta_{n}}))$, anticipating that, under appropriate conditions, these elements have an asymptotically normal distribution, with a mean of zero, in the absence of misspecification. Given a consistent estimator for the asymptotic covariance matrix, we can obtain a test statistic asymptotically $\chi^2_q$, for a specified $q$.

We now define other auxiliary matrices necessary for the construction of the misspecification test statistic. Let's consider $l=1, \ldots, p(p+1)/2; i=1, \ldots, p; j=1, \ldots, p $, where $p$ is the number of vector coordinates $\theta$ (number of model parameters). And be

 $$ d_l(U_t,\theta)=\frac{\partial \log f(U_t, \theta)}{\partial \theta_i }. \frac{\partial \log f(U_t, \theta)}{\partial \theta_j }+ \frac{\partial^2 \log f(U_t, \theta)}{\partial \theta_i\partial \theta_j}, $$
We also define
$$D_{ln}(\hat{\theta_n})=n^{-1}\sum_{t=1}^{n}d_l(U_t, \hat{\theta_n}),$$ which are the elements of $A_n(\hat{\theta_{n}})+B_n(\hat{\theta_{n}})$.

$$ D_n(\hat{\theta_{n}})=n^{-1}\sum_{t=1}^{n}d(U_t, \theta)$$


When partial derivatives and expectations exist, we define:
$$\nabla D_n(\theta)=n^{-1}\sum_{t=1}^{n}\left[\frac{\partial d_l(U_t, \theta) }{ \partial \theta_k}\right]$$
$$\nabla D(\theta)=E\left(\frac{\partial d_l(U_t, \theta)}{ \partial \theta_k}\right)$$

The following assumptions are necessary to construct a quantity with an asymptotic distribution that will be used in the construction of a hypothesis test to verify if there is a misspecification in relation of the family $H$.

\begin{supo}$\partial d_l(u, \theta) / \partial \theta_k$, $l=1, \ldots, q$, $k=1, \ldots, p$, exist and are continuous functions of $\theta$ for each $u$.
\label{supo8}
\end{supo}

\begin{supo}$|d_l(u, \theta)d_m(u, \theta) |$, $|\partial d_l(U_t, \theta) / \partial \theta_k|$, and $|d_l(u, \theta) .\partial \log(u, \theta) / \partial \theta_k| $ , $l,m=1, \ldots, q$, $k=1, \ldots, p$ are dominated by integrable functions with respect to $H$ for all $u$ and $\theta$ in $\Theta $.
\label{supo9}
\end{supo}

These assumptions play analogous to the \ref{supo4} and \ref{supo5} assumptions. The hypothesis \ref{supo8} requires continuous third derivatives for the quasi-log-likelihood function. Among other things, the hypothesis \ref{supo9} guarantees that $\nabla D(\theta)$ is finite for all $\theta$ in $\Theta$. We define:
$$V(\theta)=E\left(\left[d(U_t,\theta)-\nabla D(\theta)A(\theta)^{-1}\nabla \log f(U_t,\theta )\right].\left[d(U_t,\theta)-\nabla D(\theta)A(\theta)^{-1}\nabla \log f(U_t,\theta)\right]'\right )$$
$V(\theta_{*})$ is the asymptotic covariance matrix of $\sqrt{n}D_n(\hat{\theta_{n}})$ and we have:
\begin{supo}$V(\theta_{*})$ is non-singular.
\label{supo10}
\end{supo}
A consistent estimator for $V(\theta_{*})$ is
$$V_n(\hat{\theta_{n}})=n^{-1}\sum_{t=1}^{n}\left[d(U_t,\hat{\theta_n})-\nabla D_n (\hat{\theta_n})A_n(\hat{\theta_n})^{-1}\nabla \log f(U_t,\hat{\theta_n})\right]. \\ \left[d(U_t,\hat{\theta_n})-\nabla D_n(\hat{\theta_n})A_n(\hat{\theta_n})^{-1}\nabla \log f(U_t, \hat{\theta_n})\right]'$$

Then we have 
\begin{teo}\textbf{(Test for misspecification)} If assumptions from \ref{supo1} to \ref{supo10} are satisfied, and $g(U)=f(U, \theta_0)$, for $ \theta_0 \in \Theta$, so
\begin{itemize}
\item[i)]$\sqrt{n}D_n(\hat{\theta_{n}}) \overset{A}{\sim}N(0, V(\theta_0));$
\item[ii)]$V_n(\hat{\theta_{n}})\overset{qc}{\rightarrow}V(\theta_0)$, and $V_n(\hat{\theta_{n}})$ is almost certainly non-singular for all sufficiently large $n$;
\item[iii)]The test for misspecification:
\begin{equation}
\label{tm}\mathcal{I}_n=nD_n(\hat{\theta_{n}})'(V_n(\hat{\theta_{n}}))^{-1}D_n(\hat{\theta_{n}})
\end{equation}
has asymptotic distribution $\chi^2_q$.
\end{itemize}
\label{teo5}
\end{teo}

To perform the test, $\mathcal{I}_n$ is calculated and compared with the critical value of the distribution $\chi^2_q$ for a given test size. If (\ref{tm}) does not exceed this value, the null hypothesis that the model was specified correctly cannot be rejected.

\section{Misspecification Tests on Random Graph Models}

In this work we build misspecification tests for the ERG and SBM models from the tests developed by White (1982) presented in the previous section. These tests are important because both models have been widely used in practice to model social networks and this tool we propose can be used to verify if the model is really suitable for the database.

We will first show that the ERG and SBM models satisfy the regularity conditions described in assumptions 1 to 10. In this way, we will prove that the statistics ${\cal I}_n$ can be used to construct hypothesis tests to verify the adequacy of these models to network databases. Next, we will find the auxiliary matrices necessary for the construction of theses statistics in each case and we will present the tests for each model.

The construction of misspecification tests is possible due to two reasons: {\it First} We observe that the likelihood functions of these models can be written as a function of the probability distributions of their vertices and edges, in this way a single sample of a network can be taken as the sample of $n$ vertices or $m$ edges of this network. {\it Second} It is possible to obtain asymptotically consistent maximum likelihood estimators for the ERG and the SBM, as described in Section 2.2.

\subsection{Misspecification Tests on Exponential Random Graph Models}
In this section we will build a misspecification test for the ERG. We saw that it was possible to write the probability distribution of an ERG as a function of the probability distributions of its edges, which made it possible to obtain a quasi-likelihood function for the model. We verified that all the regularity conditions of the ERG function are valid, so it was possible to obtain a measurable quasi-likelihood estimator, asymptotically consistent estimators for the auxiliary matrices and a misspecification test for the model.

We consider the model with only one parameter, $\theta$, and with the function $T(G) = e(G)=k$ equal to the number of edges of the graph.

\begin{prop}
For the ERG the probability distribution of a graph $G$ is equal to the product of the probability distributions of its edges.
\end{prop}
We calculate the normalization constant $\varphi(\theta)$ by solving for $1 = \sum_{G \in \mathcal{G}}P( \theta ; G)$:
$$1 = \sum^{\binom{n}{2}}_{i=0}\sum_{G \in \mathcal{G}(e(G)=i)}\exp(i\theta- \ varphi(\theta))=e^{-\varphi(\theta)}\sum^{\binom{n}{2}}_{i=0}\binom{\binom{n}{2}}{ i}e^{\theta i }=e^{-\varphi(\theta)}\left(1+e^{\theta}\right)^{\binom{n}{2}}.$$
It follows that $e^{-\varphi(\theta)}=z^{-1}\left(\theta\right)=\left(1+e^{\theta}\right)^{-\binom{ n}{2}}$ . Then, the probability function of a particular graph $G$, with $k$ edges is given by:
$$ P( G; \theta )=e^{k\theta}\left(1+e^{\theta}\right)^{-\binom{n}{2}}=\left(\frac{ e^\theta}{1+e^\theta}\right)^k\left(1-\frac{e^\theta}{1+e^\theta}\right)^{\binom{n}{ 2}-k}$$

which is the expression \ref{er} of the Erdos-Rényi model with parameter $p=\left(\frac{e^\theta}{1+e^\theta}\right)$. Thus, ER models with $p\neq 0.1$ are exponential random graphs.

In this case, the probability function of a particular graph $G$, with $k$ edges is given by:
$$P(k;\theta)=e^{k\theta}(1+e^{\theta})^{-\binom{n}{2}},$$
Let's sample $U_t$, $t=1, \ldots, \binom{n}{2}$, the $\binom{n}{2}$ elements below the main diagonal of the matrix $Y$ where $ Y_{ij}$ is the random variable of $\displaystyle Bernoulli\left( \frac{e^{\theta}}{1+e^{\theta}}\right)$, variable indicating the existence of the edge between the vertices $i$ and $j$, for $i=1,\ldots,n$ and $j=1,\ldots,n$. We will take the elements below the main diagonal of $Y$ in the following order: column by column, from column $1$ to column $n-1$, from the lowest index row to the $n$ row. In this way:

\begin{small}
$$\begin{matrix}

U_1=y_{2,1} & & & & & \\
U_2=y_{3,1} & U_{(n-1)+1}=y_{3,2} & & & & \\
U_3=y_{4,1} & U_{(n-1)+2}=y_{4,2} & U_{(n-1)+(n-2)+1}=y_{4,3} & & & \\
U_4=y_{5,1} & U_{(n-1)+3}=y_{5,2} & U_{(n-1)+(n-2)+2}=y_{5,3} & U_{(n-1)+(n-2)+(n-3)+1}=y_{5,4} & & \\
U_5=y_{6,1} & U_{(n-1)+4}=y_{6,2} & U_{(n-1)+(n-2)+3}=y_{6,3} & U_{(n-1)+(n-2)+(n-3)+2}=y_{6,4} & & \\
\vdots & \vdots & \vdots & \vdots & & \\
U_{n-1}=y_{n,1} & U_{(n-1)+(n-2)}=y_{n,2} & U_{(n-1)+(n-2)+ (n-3)}=y_{n,3} & U_{(n-1)+(n-2)+(n-3)+(n-4)}=y_{n,4} & \ldots & U_{\binom{n}{2}}=y_{n, n-1}

\end{matrix}$$
\end{small}

So we have $k=\displaystyle \sum_{t=1}^{\binom{n}{2}}U_t$ and
$$P(k; \theta)=e^{k\theta}(1+e^{\theta})^{-\binom{n}{2}}=\prod_{t=1}^{\binom{n}{2}}e^{U_{t}\theta}(1+e^{\theta})^{-1}$$
that is, $\displaystyle P(k;\theta)=\prod_{t=1}^{\binom{n}{2}}f(U_t; \theta )$ where $f(U_t; \theta ) =e^{U_{t}\theta}(1+e^{\theta})^{-1}$.

\begin{prop}
The regularity conditions given in assumptions 1 to 10 are valid for $f(U_t; \theta )$ in the ERG.
\begin{enumerate}

\item As $H$ is unknown, a priori, we choose a family of distribution functions that may or may not contain the true structure, $H$. In the case of the ERG, $H$ has a distribution $\displaystyle Bernoulli\left( \frac{e^{\theta}}{1+e^{\theta}}\right)$ for all $U_t$, satisfying the assumption \ref{supo1}.

\item For the ERG, $f(U_t, \theta )=e^{U_{t}\theta}(1+e^{\theta})^{-1}$ which is measurable in $U_{t} $ for each $\theta \in \mathbb{R}$ and continuous on $\theta $ for each $U_{t} \in \left\{0,1\right\}$, satisfying the assumption \ref{supo2}.

\item $E(\log( h(U_t))$ exists and $|\log f(U_t, \theta)|=|\theta U_t-\log (1+e^{\theta})|$ for all $\theta \in \mathbb{R}$ is integrable with respect to $H$.

In this way, the assumption \ref{supo3} is satisfied.

\item $$\displaystyle \frac{\partial \log f(U_t, \theta)}{\partial \theta}=\left( U_t-\frac{e^{\theta}}{1+e^{\theta}}\right)$$ is measurable in $U_{t}$ for each $\theta \in \mathbb{R}$ and continuously differentiable in $\theta $ for each $U_{t} \in \left\{0.1\right\}$, satisfying the assumption \ref{supo4}.

\item We have $\left\vert\frac{\partial^2 \log f(u, \theta) }{ \partial \theta^{2}}\right\vert=\frac{e^{\theta} }{(1+e^{\theta})^{2}}$ and $\left\vert\frac{\partial \log f(u, \theta)}{ \partial \theta} . \frac{\partial \log f(u, \theta) }{\partial \theta} \right\vert=\left( U_t-\frac{e^{\theta}}{1+e^{\theta} }\right)^2$ are dominated by integrable functions with respect to $H$ for all $U_t$ in $\left\{0.1\right\}$ and $\theta$ in $\mathbb{R} $. So supposition \ref{supo5} satisfied.

\item In the case of the ERG $\theta_{*} \in \mathbb{R}$, $B(\theta_{*})$ is non-singular and $\theta_{*}$ is a regular point of $A( \theta)$ which satisfies the assumption \ref{supo6}.

\item Note that $\displaystyle \frac{\partial[\frac{\partial \log f(u, \theta)}{\partial \theta}. f(u,\theta)]}{\partial\theta}=\frac{e^{U_{t}\theta}(2U_t^2-2U_t-1)e^{\theta}+(U_t-1) ^2e^{2\theta}}{(1+e^{\theta})^3}$ is integrable with respect to $\nu$ for all $\theta \in \mathbb{R}$ being the assumption \ref{supo7} checked.

\item No ERG $\displaystyle \frac{\partial d_1(u, \theta) }{ \partial \theta_k} = \left(-2U_t\frac{e^{\theta}}{(1+e^{\theta})^2}+e^{\theta} \left(\frac{3e^{\theta}-1}{(1+e^{\theta})^3}\right)\right)$ is a continuous function of $\theta$ for each $U_t$, the assumption \ref{supo8} being satisfied.

\item The functions:

$$|d_1(u, \theta)d_1(u, \theta) |=\left[U_t^2-2U_t\frac{e^\theta}{1+e^\theta}+e^\theta \left (\frac{e^\theta-1}{(1+e^\theta)^2}\right)\right]^2,$$

$$\left\vert \frac{\partial d_1(U_t, \theta) }{ \partial \theta}\right\vert=\left\vert \left(-2U_t\frac{e^{\hat{\theta }}}{(1+e^{\hat{\theta}})^2}+e^{\hat{\theta}} \left(\frac{3e^{\hat{\theta}}-1 }{(1+e^{\hat{\theta}})^3}\right)\right) \right\vert \text{ and }$$

$$\left\vert d_1(u, \theta) .\frac{\partial \log(u, \theta) }{ \partial \theta}\right\vert=\left\vert \left[U_t^2- 2U_t\frac{e^\theta}{1+e^\theta}+e^\theta \left(\frac{e^\theta-1}{(1+e^\theta)^2}\right) \right]. \left[ U_t - \frac{e^{\theta}}{1+e^{\theta}} \right] \right\vert$$

are integrable with respect to $H$ for all $U_t$ and $\theta$ in $\mathbb{R}$, the assumption \ref{supo9} being satisfied.

\item In the ERG, $V(\theta_{*})$ is non-singular, and the assumption \ref{supo10} is checked.

\end{enumerate}
\end{prop}

As are valid all assumptions from \ref{supo1} to \ref{supo10}, we get a measurable QMLE, asymptotically consistent estimators for the auxiliary matrices and a misspecification test for the ERG.

\begin{teo}
Given the assumptions \ref{supo1} and \ref{supo2}, for every $n$ exists QMLE $\hat{\theta_n}$, measurable.\label{teo6}
\end{teo}

Indeed, for the ERG $\displaystyle L_n(U,\theta)\equiv \binom{n}{2}^{-1}\sum_{t=1}^{\binom{n}{2}}\ log f(U_t, \theta)=\binom{n}{2}^{-1}\sum_{t=1}^{\binom{n}{2}}\left(\theta U_t-\log ( 1+e^{\theta})\right)=$ $$=\left(\theta\frac{\sum_{t=1}^{\binom{n}{2}}U_t}{\binom{n }{2}} - log(1+e^{\theta })\right)$$ is measurable and has a maximum for $\theta \in \mathbb{R}$.

In this case the maximum is given by $\displaystyle \hat{\theta_n}=\log \left( \frac{\frac{\sum_{t=1}^{\binom{n}{2}}U_t}{\ binom{n}{2}}}{1-\frac{\sum_{t=1}^{\binom{n}{2}}U_t}{\binom{n}{2}}}\right). $

Let $\hat{\theta}=\hat{\theta_n}$, let's define the following auxiliary matrices:

$$A_n(\hat{\theta})=\binom{n}{2}^{-1}\sum_{t=1}^{\binom{n}{2}}\left[\frac{\partial^2 \log f(U_t, \hat{\theta})}{\partial \theta^2}\right]=\binom{n}{2}^{-1}\sum_{t=1}^ {\binom{n}{2}}\left( -\frac{e^{\hat{\theta}}}{(1+e^{\hat{\theta}})^{2}}\right )=-\frac{e^{\hat{\theta}}}{(1+e^{\hat{\theta}})^{2}},$$

$$B_n(\hat{\theta})=\binom{n}{2}^{-1}\sum_{t=1}^{\binom{n}{2}}\left[\frac{\partial \log f(U_t, \hat{\theta})}{\partial \theta}. \frac{\partial \log f(U_t, \hat{\theta})}{\partial \theta}\right]=\binom{n}{2}^{-1}\sum_{t=1}^ {\binom{n}{2}}\left( \left( U_t-\frac{e^{\hat{\theta}}}{1+e^{\hat{\theta}}}\right)^ 2\right).$$

$$C_n(\hat{\theta})=A_n(\hat{\theta})^{-1}B_n(\hat{\theta})A_n(\hat{\theta})^{-1}= \frac{\left(1+e^{\hat{\theta}}\right)^4}{e^{2\hat{\theta}}}\binom{n}{2}^{-1} \sum_{t=1}^{\binom{n}{2}}\left( \left( U_t-\frac{e^{\hat{\theta}}}{1+e^{\hat{\theta}}}\right)^2\right),$$

Let's define
$$d_l(U,\theta)=\frac{\partial \log f(U_t, \theta)}{\partial \theta_i }. \frac{\partial \log f(U_t, \theta)}{\partial \theta_j }+ \frac{\partial^2 \log f(U_t, \theta)}{\partial \theta_i\partial \theta_j}, $$
$l=1, \ldots, p(p+1)/2; i=1, \ldots, p; j=1, \ldots, p $. Where $p$ is the number of coordinates of the vector $\theta$.

In the case of the ERG, we only have one parameter, so we calculate the $d_{1}$, given by

$$d_1(U_t,\theta)=\frac{\partial \log f(U_t, \theta)}{\partial \theta }. \frac{\partial \log f(U_t, \theta)}{\partial \theta }+ \frac{\partial^2 \log f(U_t, \theta)}{\partial \theta^2}=U_t^ 2-2U_t\frac{e^\theta}{1+e^\theta}+e^\theta \left(\frac{e^\theta-1}{(1+e^\theta)^2}\right).$$

The test will be based on $\displaystyle D_{n}(\hat{\theta})=\binom{n}{2}^{-1}\sum_{t=1}^{\binom{n}{2 }}d_l(U_t, \hat{\theta})$, which are the elements of $A_n(\hat{\theta})+B(\hat{\theta})$.

Let $l=1, \ldots\, q=p(p+1)/2$, define the vector $d(U_t,\theta)$, of dimension $q \times 1$, like this

$$d(U_t,\theta)=U_t^2-2U_t\frac{e^\theta}{1+e^\theta}+e^\theta \left(\frac{e^\theta-1}{ (1+e^\theta)^2}\right)$$

So $D_n(\hat{\theta})=\binom{n}{2}^{-1}\sum_{t=1}^{\binom{n}{2}}d(U_t, \hat{\theta})$. For the ERG,

$$D_n(\hat{\theta})=\binom{n}{2}^{-1}\sum_{t=1}^{\binom{n}{2}}\left(U_t^2- 2U_t\frac{e^{\hat{\theta}}}{1+e^{\hat{\theta}}}+e^{\hat{\theta}} \left(\frac{e^{\hat{\theta}}-1}{(1+e^{\hat{\theta}})^2}\right)\right)$$

 From $D_n(\hat{\theta})$, we define:
$$\nabla D_n(\theta)=\binom{n}{2}^{-1}\sum_{t=1}^{\binom{n}{2}}\left[\frac{\partial d_1 (U_t, \hat{\theta}) }{ \partial \theta}\right]=\binom{n}{2}^{-1}\sum_{t=1}^{\binom{n}{2 }}\left(-2U_t\frac{e^{\hat{\theta}}}{(1+e^{\hat{\theta}})^2}+e^{\hat{\theta}} \left(\frac{3e^{\hat{\theta}}-1}{(1+e^{\hat{\theta}})^3}\right)\right)$$

$V(\theta_{*})$ is the asymptotic covariance matrix of $\sqrt{n}D_n(\hat{\theta})$ and we have that a consistent estimator for $V(\theta_{*})$ is
{\small $$V_n(\hat{\theta})=\binom{n}{2}^{-1}\sum_{t=1}^{\binom{n}{2}}\left[d (U_t,\hat{\theta})-\nabla D_n(\hat{\theta})A_n(\hat{\theta})^{-1}\nabla \log f(U_t,\hat{\theta} )\right]. \\ \left[d(U_t,\hat{\theta})-\nabla D_n(\hat{\theta})A_n(\hat{\theta})^{-1}\nabla \log f(U_t, \hat{\theta})\right]'.$$}

In the ERG:

$$V_n(\hat{\theta})=\binom{n}{2}^{-1}\sum_{t=1}^{\binom{n}{2}}\left[\left(U_t ^2-2U_t\frac{e^{\hat{\theta}}}{1+e^{\hat{\theta}}}+e^{\hat{\theta}} \left(\frac{e ^{\hat{\theta}}-1}{(1+e^{\hat{\theta}})^2}\right)\right) \right. +$$ $$ \left. + \left(\binom{n}{2}^{-1}\sum_{t=1}^{\binom{n}{2}}\left(U_t^2-2U_t\frac{e^{\ hat{\theta}}}{1+e^{\hat{\theta}}}+e^{\hat{\theta}} \left(\frac{e^{\hat{\theta}}-1 }{(1+e^{\hat{\theta}})^2}\right)\right)\right). \frac{(1+e^{\hat{\theta}})^{2}}{e^{\hat{\theta}}}.\left( U_t-\frac{e^{\hat{\ theta}}}{1+e^{\hat{\theta}}}\right)\right]^{2}$$

Thus, we generate the following test:

\begin{teo}\textbf{(Misspecification Test on Exponential Random Graphs)} Assumptions from \ref{supo1} to \ref{supo10} are satisfied, and if $h(U)=f(U, \theta_0) $, to $\theta_0 \in \Theta$, then
\begin{equation}
\label{tmge}\mathcal{I}_n=\frac{1}{V_n(\hat{\theta})}\left(\binom{n}{2}^{-1}\sum_{t=1}^{\binom{n}{2}}\left(U_t^2-2U_t\frac{e^{\hat{\theta}}}{1+e^{\hat{\theta}}}+e ^{\hat{\theta}} \left(\frac{e^{\hat{\theta}}-1}{(1+e^{\hat{\theta}})^2}\right)\right)\right)^2
\end{equation}
has asymptotic distribution $\chi^2_1$.
\label{teo7}
\end{teo}

For a hypothesis test with $\alpha$ of significance, calculate \ref{tmge} and use the criterion: if $\mathcal{I}_n \leq \chi^2_{(\alpha, 1) } $, where $\chi^2_{(\alpha, 1)}$ is the cumulative distribution value of $\chi^2_{1}$ in $\alpha$, the model was well specified, otherwise the model was misspecified.

\subsection{Misspecification Test on Stochachistic Block Models}
In this section we will build a specification error test for SBM. We saw that it was possible to write the probability distribution of an SBM as a function of the probability distributions of its edges, when the classes and the number of connections of each of its vertices are given, which made it possible to obtain a quasi- likelihood for the model. We verified that all regularity conditions of the SBM function are valid, so it was possible to obtain a measurable quasi-likelihood estimator, asymptotically consistent estimators for the auxiliary matrices and a specification error test for the SBM.

The probability distribution function for a graph in SBM is given by:
$$P( \theta, \eta ; X,Y)=\theta_{1}^{n_1} \cdots \theta_{m}^{n_m}\prod_{1 \leq k \leq l \leq m}\eta_{kl}^{e_{kl}}(1-\eta_{kl})^{n_{kl}-e_{kl}},$$
where $n_k=\sum_{i=1}^{n}I(X_i=k)$ denotes the number of vertices of $G$ that belong to block $k$,
$$e_{kl}=\sum_{1 \leq i \neq j \leq n}Y_{ij}I(x_i=k)I(x_j=l).$$
denotes the number of edges of $G$ that have a vertex in block $k$ and a vertex in block $j$, and
$$n_{kl}=\left\{ \begin{array}{cc}
n_k n_l & \text{ if } k \neq l \\
\binom{n_k}{2} & \text{ if } k=l,
\end{array}\right.$$

To build the Misspecification Test in the case of SBM we need to rewrite $P( \theta, \eta ; X,Y)$ as a function of the random variable $U_t$. The following proposition shows how to do this.

\begin{prop}
For SBM, the probability distribution of a graph $G$ is equal to the product of the probability distributions of its edges, when given the blocks and the number of connections at each of its vertices. So,

  $$P( \theta, \eta ; X,Y)= f(U_t,\theta, \eta )=\left(\theta_k^{\frac{\textbf{I}_{x_{i}=k} }{ n_i }} \theta_l^{\frac{\textbf{I}_{x_{j}=l}}{ n_j }}\right) \eta_{kl}^{y_{ij}}(1-\eta_{kl})^{ \left( 1 -y_{ij} \right)}.$$
\end{prop}

Let's construct the sample vector $U_t$, $t=1, \ldots, \binom{n}{2}$, for each edge of the graph $G$, we will associate to each edge $t$ two vertices $i^t $(vertex $i$ of edge $t$) and $j^t$(vertex $j$ of edge $t$) , $U_t$ will be obtained from the sample matrices $X$ and $Y$.

The vector $U_t$, of $5$ dimensions, for each edge $t$, must contain the following information:
\begin{itemize}
\item[$\bullet$]the block $k$ of the vertex $i^t$, we will denote this block by $k^{t}$;
\item[$\bullet$]the block $l$ of the vertex $j^t$, we will denote this block by $l^{t}$;
\item[$\bullet$]the number of links from vertex $i^t$: $n_i^t=\sum_{p=1}^{n}y_{ip}$;
\item[$\bullet$]the number of links from the vertex $j^t$: $n_j^t=\sum_{p=1}^{n}y_{jp}$;
\item[$\bullet$]the random variable $y_{ij}^t$ of $Bernoulli( \eta_{kl})$( variable indicating the edge $t$ between the vertices $i^t$ and $j^ t$.)
\end{itemize}

That is, $U_t=\left( k^t,l^t,n_i^t , n_j^t, y_{ij}^t\right)$, $t=1, \ldots, \binom{n}{ 2}$.

This construction of $U_t$ is plausible because is possible associate to each edge $t$ two vertices, $i^t$(vertex $i$ of edge $t$) and $j^t$(vertex $j$ of edge $t$), just notice that for the $\binom{n}{2}$ elements $y_{ij}$ that are below the main diagonal of the matrix $Y$, their vertex $i^t$ is in the block $x_i$( where $i$ is the row index of the element $y_{ij}$) and the vertex $j^t$ is in the block $x_j$ ( where $j$ is the column index of the element $y_ {ij}$), for each $t$. Thus, once again, we will take the elements below the main diagonal of $Y$ in the following order: column by column, from column $1$ to column $n-1$, in the direction of the row with the lowest index to the row $n$ . So:

$$\begin{matrix}

U_1=(x_2, x_1, n_2, n_1, y_{2,1}) & & \\
U_2=(x_3, x_1, n_3, n_1, y_{3,1}) & U_{(n-1)+1}=(x_3, x_2, n_3, n_2, y_{3,2}) & \\
U_3=(x_4, x_1, n_4, n_1, y_{4,1}) & U_{(n-1)+2}=(x_4, x_2, n_4, n_2, y_{4,2}) & U_{(n-1)+(n-2)+1}=(x_4, x_3, n_4, n_3, y_{4,3}) \\
U_4=(x_5, x_1, n_5, n_1, y_{5,1}) & U_{(n-1)+3}=(x_5, x_2, n_5, n_2, y_{5,2}) & U_{(n-1)+(n-2)+2}=(x_5, x_3, n_5, n_3, y_{5,3}) \\
U_5=(x_6, x_1, n_6, n_1, y_{6,1}) & U_{(n-1)+4}=(x_6, x_2, n_6, n_2, y_{6,2}) & U_{(n-1)+(n-2)+3}=(x_6, x_3, n_6, n_3, y_{6,3}) \\
\vdots & \vdots & \vdots \\
U_{n-1}=(x_n, x_1, n_n, n_1, y_{n,1}) & U_{(n-1)+(n-2)}=(x_n, x_2, n_n, n_2, y_{n,2}) & U_{(n-1)+(n-2)+(n-3)}=(x_n, x_3, n_n, n_3, y_{n,3}) \\

\end{matrix}$$

$$\begin{matrix}
& & \\
U_{(n-1)+(n-2)+(n-3)+1}=(x_5, x_4, n_5, n_4, y_{5,4}) & & \\
 U_{(n-1)+(n-2)+(n-3)+2}=(x_6, x_4, n_6, n_4, y_{6,4}) &  &  \\
 \vdots & & \\
 U_{(n-1)+(n-2)+(n-3)+(n-4)}=(x_n, x_4, n_n, n_4, y_{n,4}) & \ldots & U_{\binom{n}{2}}=(x_n, x_{(n-1)}, n_n, n_{(n-1)}, y_{n,n-1})
\end{matrix}$$

Hence, we have 
$$P( \theta, \eta ; X,Y)=\theta_{1}^{n_1} \cdots \theta_{m}^{n_m}\prod_{1 \leq k \leq l \leq m}\eta_{kl}^{e_{kl}}(1-\eta_{kl})^{n_{kl}-e_{kl}}=\prod_{t=1}^{\binom{n}{2}}
\theta_k^{\frac{\textbf{I}_{x_{i^t}=k}}{ n_i^t  }} \theta_l^{\frac{\textbf{I}_{x_{j^t}=l}}{ n_j^t  }} \eta_{kl}^{y_{ij}^t}(1-\eta_{kl})^{ \left( 1 -y_{ij}^t \right)}
$$

that is, $\displaystyle P( \theta, \eta ; X,Y)=\prod_{t=1}^{\binom{n}{2}}f(U_t; \theta, \eta )$ where $f(U_t,\theta, \eta )=\left(\theta_k^{\frac{\textbf{I}_{x_{i^t}=k}}{ n_i^t  }} \theta_l^{\frac{\textbf{I}_{x_{j^t}=l}}{ n_j^t  }}\right) \eta_{kl}^{y_{ij}^t}(1-\eta_{kl})^{ \left( 1 -y_{ij}^t \right)}$

To simplify the notation a bit let's use $U_t=\left( k^t,l^t,n_i^t , n_j^t, y_{ij}^t\right)=\left( k,l ,n_i , n_j , y_{ij} \right)$, that is, let's omit the subscript $t$ in all elements of $U_t$, knowing that each one of them depends on $t$.

So $f(U_t,\theta, \eta )=\left(\theta_k^{\frac{\textbf{I}_{x_{i}=k}}{ n_i }} \theta_l^{\frac{ \textbf{I}_{x_{j}=l}}{ n_j }}\right) \eta_{kl}^{y_{ij}}(1-\eta_{kl})^{ \left( 1 - y_{ij} \right)}$

\begin{prop}
The regularity conditions given in assumptions 1 to 10 are valid for $f(U_t,\theta, \eta )$ in SBM.
\begin{enumerate}

\item As $H$ is unknown a priori, we choose a family of distribution functions that may or may not contain the true structure, $H$. For SBM, $H$ has joint distribution $P( \theta, \eta ; X,Y)$ with $P(x_i=k)=\theta_k \in (0,1)$ and $ \eta_{kl} \in (0,1)$, for all $k,l=1, \ldots, m$, satisfying the assumption \ref{supo1}.

To simplify the notation, we will denote the set $(\theta, \eta)$ just by $\theta$ from this point forward.

\item For the ERG, $f(U_t, \theta )=\left( \theta_k^{\frac{\textbf{I}_{x_i=k}}{n_i} } \theta_l^{\frac{\textbf {I}_{x_j=l}}{n_j} } \right) \eta_{kl}^{y_{ij}}(1-\eta_{kl})^{\left( 1 -y_{ij} \right)}$ is measurable in $U_{t}$ for each $\theta \in \Theta$ and continuous in $\theta $ for each $U_{t} \in \Omega $, satisfying the assumption \ref{supo2}.

\item $E(\log( h(U_t))$ exists and $$|\log f(U_t, \theta)|=\left\vert\left(\frac{\textbf{I}_{x_i=k }}{n_i} \right)\log (\theta_k) +\left(\frac{\textbf{I}_{x_j=l}}{n_j} \right)\log (\theta_l)+y_{ij} \log (\eta_{kl})+ \left( 1 -y_{ij} \right)\log (1-\eta_{kl}))\right\vert$$ for all $\theta \in \Theta$ is integrable with respect to $H$. In this way, the assumption \ref{supo3} is satisfied.

\item The functions
$$\frac{\partial \log f(U_t,\theta )}{\partial \theta_k} =\frac{\left(\frac{\textbf{I}_{x_i=k}}{n_i} \right )}{\theta_k}$$
$$\frac{\partial \log f(U_t,\theta )}{\partial \theta_l} =\frac{\left(\frac{\textbf{I}_{x_j=l}}{n_j}\right )}{\theta_l}$$
$$\frac{\partial \log f(U_t, \theta )}{\partial \eta_{kl}} =\frac{y_{ij}}{\eta_{kl}}-\frac{\left( 1 -y_{ij} \right)}{1-\eta_{kl}}$$

are measurable in $U_{t}$ for each $\theta \in \Theta$ and continuously differentiable in $\theta $ for each $U_{t} \in \Omega$, satisfying the assumption \ref{supo4}.

\item Let's see that the modulus of the products of the first-order derivatives and the modulus of the second-order derivatives of $\log f(U_t,\theta )$ with respect to each of the parameters, are dominated by integrable functions with respect to $ H$ for all $U_t$ in $\Omega$ and $\theta$ in $\Theta$. So the assumption \ref{supo5} will be satisfied. Following are the functions:

$$\left\vert \left(\frac{\partial \log f(U_t, \theta )}{\partial \theta_k}\right)^2 \right\vert  = \left\vert \frac{\left(\frac{\textbf{I}_{x_i=k}}{n_i}\right)^2}{\theta_k^2} \right\vert ,
 \left\vert\frac{\partial \log f(U_t, \theta )}{\partial \theta_k}\frac{\partial \log f(U_t, \theta )}{\partial \theta_l} \right\vert  =\left\vert\frac{\left(\frac{\textbf{I}_{x_i=k}}{n_i}\right)\left(\frac{\textbf{I}_{x_j=l}}{n_j}\right)}{\theta_k\theta_l}\right\vert$$
{\small $$\left\vert\frac{\partial \log f(U_t, \theta )}{\partial \theta_k} \frac{\partial \log f(U_t, \theta )}{\partial \eta_{kl}}\right\vert  =\left\vert\frac{\left(\frac{\textbf{I}_{x_i=k}}{n_i}\right)}{\theta_k}\left( \frac{y_{ij}}{\eta_{kl}}-\frac{\left(1 -y_{ij} \right)}{1-\eta_{kl}}\right)\right\vert ,\left\vert \frac{\partial \log f(U_t, \theta )}{\partial \theta_l}\frac{\partial \log f(U_t, \theta )}{\partial \theta_k} \right\vert =  \left\vert \frac{\left(\frac{\textbf{I}_{x_i=k}}{n_i}\right)\left(\frac{\textbf{I}_{x_j=l}}{n_j}\right)}{\theta_l\theta_k}\right\vert$$} 

$$\left\vert\left(\frac{\partial \log f(U_t, \theta )}{\partial \theta_l}\right)^2  \right\vert=\left\vert\frac{\left(\frac{\textbf{I}_{x_j=l}}{n_j}\right)^2}{\theta_l^2}\right\vert , \left\vert\frac{\partial \log f(U_t, \theta )}{\partial \theta_l} \frac{\partial \log f(U_t, \theta )}{\partial \eta_{kl}} \right\vert =\left\vert\frac{\left(\frac{\textbf{I}_{x_j=l}}{n_j}\right)}{\theta_l}\left( \frac{y_{ij}}{\eta_{kl}}-\frac{\left( 1 -y_{ij} \right)}{1-\eta_{kl}}\right)\right\vert$$

$$\left\vert\frac{\partial \log f(U_t, \theta )}{\partial \eta_{kl}}\frac{\partial \log f(U_t, \theta )}{\partial \theta_k}\right\vert   =\left\vert\left( \frac{y_{ij}}{\eta_{kl}}-\frac{\left( 1 -y_{ij} \right)}{1-\eta_{kl}}\right)\frac{\left(\frac{\textbf{I}_{x_i=k}}{n_i}\right)}{\theta_k}\right\vert$$ 

$$\left\vert\frac{\partial \log f(U_t, \theta )}{\partial \eta_{kl}} \frac{\partial \log f(U_t, \theta )}{\partial \theta_l} \right\vert= \left\vert\left( \frac{y_{ij}}{\eta_{kl}}-\frac{\left( 1 -y_{ij} \right)}{1-\eta_{kl}}\right)\frac{\left(\frac{\textbf{I}_{x_j=l}}{n_j}\right)}{\theta_l}\right\vert , \left\vert\left(\frac{\partial \log f(U_t, \theta )}{\partial \eta_{kl}}\right)^2 \right\vert =\left\vert\left(\frac{y_{ij}}{\eta_{kl}}-\frac{\left(1-y_{ij} \right)}{1-\eta_{kl}} \right)^2\right\vert$$

$$\left\vert\frac{\partial^2 \log f(U_t, \theta )}{\partial \theta_k^2}\right\vert =\left\vert\frac{-\left(\frac{\textbf{I}_{x_i=k}}{n_i}\right)}{\theta_k^2}\right\vert , \left\vert\frac{\partial^2 \log f(U_t, \theta )}{\partial \theta_k\partial \theta_l}\right\vert =0 ,\left\vert 
\frac{\partial^2 \log f(U_t, \theta )}{\partial \theta_k\partial \eta_{kl}} \right\vert=0$$
$$\left\vert\frac{\partial^2 \log f(U_t, \theta )}{\partial \theta_l\partial \theta_k}\right\vert =0 , \left\vert\frac{\partial^2 \log f(U_t, \theta )}{\partial \theta_l^2}\right\vert =\left\vert \frac{-\left(\frac{\textbf{I}_{x_j=l}}{n_j}\right)}{\theta_l^2}\right\vert , \left\vert\frac{\partial^2 \log f(U_t, \theta )}{\partial \theta_l\partial \eta_{kl}}\right\vert =0$$

$$\left\vert\frac{\partial^2 \log f(U_t, \theta )}{\partial \eta_{kl}\partial \theta_k}\right\vert =0, \left\vert\frac{\partial^2 \log f(U_t, \theta )}{\partial \eta_{kl}\partial \theta_l}\right\vert =0 , \left\vert\frac{\partial^2 \log f(U_t, \theta )}{\partial \eta_{kl}^2}\right\vert =\left\vert\frac{-y_{ij}}{\eta_{kl}^2}-\frac{ \left( 1 -y_{ij} \right)}{(1-\eta_{kl})^2}\right\vert$$

which are all dominated by integrable functions with respect to $H$.

\item In the case of SBM $\theta_{*} \in \mathbb{R}$, $B(\theta_{*})$ is non-singular and $\theta_{*}$ is a regular point of $A( \theta)$ which satisfies the \ref{supo6} assumption.

\item Note that

$$\displaystyle \frac{\partial[\frac{\partial \log f(U_t, \theta)}{\partial \theta_k}. f(U_t,\theta)]}{\partial\theta_k}= \left(\frac{\textbf{I}_{x_i=k}}{n_i} \right)\left(\frac{\textbf{I}_{x_i=k}}{n_i} - 1 \right) \left( \theta_k^{\left(\frac{\textbf{I}_{x_i=k}}{n_i} -2 \right) } \theta_l^{\frac{\textbf{I}_{x_j=l}}{n_j} }\right) \eta_{kl}^{y_{ij}}(1-\eta_{kl})^{\left( 1 -y_{ij} \right)}$$ 

$$\displaystyle \frac{\partial[\frac{\partial \log f(U_t, \theta)}{\partial \theta_l}. f(U_t,\theta)]}{\partial\theta_l}=  \left(\frac{\textbf{I}_{x_j=l}}{n_j} \right)\left(\frac{\textbf{I}_{x_j=l}}{n_j} - 1 \right) \left( \theta_k^{\frac{\textbf{I}_{x_i=k}}{n_i} }\theta_l^{\left(\frac{\textbf{I}_{x_j=l}}{n_i} -2 \right) } \right) \eta_{kl}^{y_{ij}}(1-\eta_{kl})^{\left( 1 -y_{ij} \right)}$$

\begin{footnotesize}
$$\displaystyle \frac{\partial[\frac{\partial \log f(U_t, \theta)}{\partial \eta_{kl}}. f(U_t,\theta)]}{\partial\eta_{kl}}= \left(\theta_k^{\frac{\textbf{I}_{x_i=k}}{n_i}  } \theta_l^{\frac{\textbf{I}_{x_j=l}}{n_j} }\right)\left[  \left(\frac{\eta_{kl}^{y_{ij}-1}\left(y_{ij}-\eta_{kl}\right)}{1-\eta_{kl}}\right).\left( \frac{y_{ij}-\eta_{kl}}{\eta_{kl}-\eta_{kl}^2}\right) + \left(\frac{\left(1-\eta_{kl}\right)\eta_{kl}^{y_{ij}}}{1-\eta_{kl}} \right)\left(\frac{y_{ij}\left(2\eta_{kl}-1\right)-\eta_{kl}^2}{\left(1-\eta_{kl}\right)^2\eta_{kl}}\right) \right]$$
\end{footnotesize}

are integrable with respect to $\nu$ for all $\theta \in \Theta$ and the assumption \ref{supo7} is checked.

\item In SBM \\
$\displaystyle \frac{\partial d_1(u, \theta) }{ \partial \theta_k} = \left( \frac{-2\left(\frac{\textbf{I}_{x_i=k}}{n_i}\right)^2}{\theta_k^3}+\frac{2\left(\frac{\textbf{I}_{x_i=k}}{n_i}\right)}{\theta_k^3} \right)$, $\displaystyle \frac{\partial d_1(u, \theta) }{ \partial \theta_l} = 0 $, $\displaystyle \frac{\partial d_1(u, \theta) }{ \partial \eta_{kl}} = 0 $, 

$\displaystyle \frac{\partial d_2(u, \theta) }{ \partial \theta_k} = \left( \frac{-\left(\frac{\textbf{I}_{x_i=k}}{n_i}\right)\left(\frac{\textbf{I}_{x_j=l}}{n_j}\right)}{\theta_k^2\theta_l} \right)$, $\displaystyle \frac{\partial d_2(u, \theta) }{ \partial \theta_l} =  \left( \frac{-\left(\frac{\textbf{I}_{x_i=k}}{n_i}\right)\left(\frac{\textbf{I}_{x_j=l}}{n_j}\right)}{\theta_k\theta_l^2} \right)$, $\displaystyle \frac{\partial d_2(u, \theta) }{ \partial \eta_{kl}} = 0$,  

$\displaystyle \frac{\partial d_3(u, \theta) }{ \partial \theta_k} =  \left( \frac{-\left(\frac{\textbf{I}_{x_i=k}}{n_i}\right)}{\theta_k^2} \left( \frac{y_{ij}}{\eta_{kl}}-\frac{\left( 1 -y_{ij} \right)}{1-\eta_{kl}}\right) \right)$,
$\displaystyle \frac{\partial d_3(u, \theta) }{ \partial \theta_l} = 0$,
$\displaystyle \frac{\partial d_3(u, \theta) }{ \partial \eta_{kl}} = \left(\frac{\left(\frac{\textbf{I}_{x_i=k}}{n_i}\right)}{\theta_k}\left( \frac{-y_{ij}}{\eta_{kl}^2}-\frac{\left( 1 -y_{ij} \right)}{(1-\eta_{kl})^2}\right) \right)$,

$\displaystyle \frac{\partial d_4(u, \theta) }{ \partial \theta_k} =0 $,
$\displaystyle \frac{\partial d_4(u, \theta) }{ \partial \theta_l} =\left(\frac{-2\left(\frac{\textbf{I}_{x_j=l}}{n_j}\right)^2}{\theta_l^3} +  \frac{2\left(\frac{\textbf{I}_{x_j=l}}{n_j}\right)}{\theta_l^3} \right)  $,
$\displaystyle \frac{\partial d_4(u, \theta) }{ \partial \eta_{kl}} =0 $,

$\displaystyle \frac{\partial d_5(u, \theta) }{ \partial \theta_k} =0 $,
$\displaystyle \frac{\partial d_5(u, \theta) }{ \partial \theta_l} = \left(-\frac{\left(\frac{\textbf{I}_{x_j=l}}{n_j}\right)}{\theta_l^2}\left( \frac{y_{ij}}{\eta_{kl}}-\frac{\left( 1 -y_{ij} \right)}{1-\eta_{kl}}\right)\right) $,
$\displaystyle \frac{\partial d_5(u, \theta) }{ \partial \eta_{kl}} =  \left(\frac{\left(\frac{\textbf{I}_{x_j=k}}{n_j}\right)}{\theta_l}\left( \frac{-y_{ij}}{\eta_{kl}^2}-\frac{\left( 1 -y_{ij} \right)}{(1-\eta_{kl})^2}\right) \right)$,

$\displaystyle \frac{\partial d_6(u, \theta) }{ \partial \theta_k} =0 $,
$\displaystyle \frac{\partial d_6(u, \theta) }{ \partial \theta_l} = 0$,
$\displaystyle \frac{\partial d_6(u, \theta) }{ \partial \eta_{kl}} =\left( 2\left(\frac{y_{ij}}{\eta_{kl}}-\frac{\left( 1 -y_{ij} \right)}{1-\eta_{kl}} \right) \left( \frac{-y_{ij}}{\eta_{kl}^2}-\frac{\left( 1 -y_{ij} \right)}{(1-\eta_{kl})^2}\right)   + \frac{2y_{ij}}{\eta_{kl}^3}-\frac{ 2\left( 1 -y_{ij} \right)}{(1-\eta_{kl})^3}\right) $,
 
are continuous functions of $\theta \in \Theta$ for each $U_t \in \Omega$, the assumption \ref{supo8} being satisfied.

\item All functions:

$\left\vert d_i(u, \theta)d_j(u, \theta) \right\vert$,
$\left\vert \frac{\partial d_i(U_t, \theta) }{ \partial \theta_k}\right\vert$, $\left\vert \frac{\partial d_i(U_t, \theta) }{ \partial \theta_l}\right\vert$, $\left\vert \frac{\partial d_i(U_t, \theta) }{ \partial \eta_{kl}}\right\vert$,
$\left\vert d_i(u, \theta) .\frac{\partial \log(u, \theta) }{ \partial \theta_k}\right\vert$,
$\left\vert d_i(u, \theta).\frac{\partial \log(u, \theta) }{ \partial \theta_l}\right\vert$,
$\left\vert d_i(u, \theta) .\frac{\partial \log(u, \theta) }{ \partial \eta_{kl}}\right\vert$,

are products of functions integrable with respect to $H$, so these functions are integrable with respect to $H$ for all $1\leq i\leq 6$, $U_t$ and $\theta \in \Theta$, satisfying the guess \ref{supo9}.

\item In SBM, $V(\theta_{*})$ is non-singular, and the assumption \ref{supo10} is checked.
\end{enumerate}
\end{prop}

As the assumptions from \ref{supo1} to \ref{supo10} are valid, we obtain a measurable QMLE, asymptotically consistent estimators for the auxiliary matrices and a misspecification test for the SBM.

\begin{teo}
Given the assumptions \ref{supo1} and \ref{supo2}, for every $n$ there is a measurable QMLE $\hat{\theta_n}$.\label{teo8}
\end{teo}

Indeed for SBM, $\displaystyle L_n(U,\theta)\equiv \binom{n}{2}^{-1}\sum_{t=1}^{\binom{n}{2}}\ log f(U_t, \theta)=$

$$\binom{n}{2}^{-1}\sum_{t=1}^{\binom{n}{2}} \left[\left(\frac{\textbf{I}_{x_i =k}}{n_i} \right)\log (\theta_k) +\left(\frac{\textbf{I}_{x_j=l}}{n_j} \right)\log (\theta_l)+y_{ ij}\log (\eta_{kl})+ \left( 1 -y_{ij} \right)\log (1-\eta_{kl})\right] $$
is measurable and has a maximum for $\theta \in \Theta$.

In this case the maximum can be obtained via Variational EM and in this work we will use the estimates via Variational EM, obtained by the package \textit{blockmodels} of the R software.

Let $\hat{\theta}=\hat{\theta_n}$, let's define the following auxiliary matrices:
$$\displaystyle \nabla \log f(U_t, \theta)= \left( \begin{matrix}
\frac{\partial \log f(U_t,\theta )}{\partial \theta_k} \\
\frac{\partial \log f(U_t,\theta )}{\partial \theta_l} \\
\frac{\partial \log f(U_t,\theta )}{\partial \eta_{kl}}
\end{matrix} \right)=\left( \begin{matrix} 
\frac{\left(\frac{\textbf{I}_{x_i=k}}{n_i} \right)}{\theta_k}\\
\frac{\left(\frac{\textbf{I}_{x_j=l}}{n_j}\right)}{\theta_l} \\
\left(\frac{y_{ij}}{\eta_{kl}}-\frac{\left( 1 -y_{ij} \right)}{1-\eta_{kl}}\right)
\end{matrix} \right)$$

$\displaystyle A_n(\hat{\theta})=\left(\begin{matrix}
\binom{n}{2}^{-1}\sum_{t=1}^{\binom{n}{2}}\left[\frac{\partial^2 \log f(U_t, \hat{\theta})}{\partial \theta_k^2}\right] & \binom{n}{2}^{-1}\sum_{t=1}^{\binom{n}{2}}\left[\frac{\partial^2 \log f(U_t, \hat{\theta})}{\partial \theta_k\partial \theta_l}\right] & \binom{n}{2}^{-1}\sum_{t=1}^{\binom{n}{2}}\left[\frac{\partial^2 \log f(U_t, \hat{\theta})}{\partial \theta_k\partial \eta_{kl}}\right] \\ 
\binom{n}{2}^{-1}\sum_{t=1}^{\binom{n}{2}}\left[\frac{\partial^2 \log f(U_t, \hat{\theta})}{\partial \theta_l\partial \theta_k}\right] & \binom{n}{2}^{-1}\sum_{t=1}^{\binom{n}{2}}\left[\frac{\partial^2 \log f(U_t, \hat{\theta})}{\partial \theta_l^2}\right] & \binom{n}{2}^{-1}\sum_{t=1}^{\binom{n}{2}}\left[\frac{\partial^2 \log f(U_t, \hat{\theta})}{\partial \theta_l\partial \eta_{kl}}\right] \\
\binom{n}{2}^{-1}\sum_{t=1}^{\binom{n}{2}}\left[\frac{\partial^2 \log f(U_t, \hat{\theta})}{\partial \eta_{kl}\partial \theta_k}\right] & \binom{n}{2}^{-1}\sum_{t=1}^{\binom{n}{2}}\left[\frac{\partial^2 \log f(U_t, \hat{\theta})}{\partial \eta_{kl}\partial \theta_l}\right] & \binom{n}{2}^{-1}\sum_{t=1}^{\binom{n}{2}}\left[\frac{\partial^2 \log f(U_t, \hat{\theta})}{\partial \eta_{kl}^2}\right] 

\end{matrix} \right)$
$$\displaystyle A_n(\hat{\theta})=\left(\begin{matrix}
\binom{n}{2}^{-1}\sum_{t=1}^{\binom{n}{2}}\left(  \frac{-\left(\frac{\textbf{I}_{x_i=k}}{n_i}\right)}{\hat{\theta_k}^2}\right) & 0 & 0 \\ 
0 & \binom{n}{2}^{-1}\sum_{t=1}^{\binom{n}{2}}\left( \frac{-\left(\frac{\textbf{I}_{x_j=l}}{n_j}\right)}{\hat{\theta_l}^2} \right) & 0 \\
0 & 0 & \binom{n}{2}^{-1}\sum_{t=1}^{\binom{n}{2}}\left( \frac{-y_{ij}}{\hat{\eta_{kl}}^2}-\frac{ \left( 1 -y_{ij} \right)}{(1-\hat{\eta_{kl}})^2}\right) 

\end{matrix} \right)$$

let's define
$$d_l(U,\theta)=\frac{\partial \log f(U_t, \theta)}{\partial \theta_i }. \frac{\partial \log f(U_t, \theta)}{\partial \theta_j }+ \frac{\partial^2 \log f(U_t, \theta)}{\partial \theta_i\partial \theta_j}, $$
$l=1, \ldots, p(p+1)/2; i=1, \ldots, p; j=1, \ldots, p $. Where $p$ is the number of coordinates of the vector $\theta$.

In the case of SBM, we have the parameters $\theta_k, \theta_l \text{ and } \eta_{kl}$ for each $U_t$, that is, 3 parameters, in this way we calculate the $d_{l}$, for $l =1,2, \ldots , 6$ given by

$$d_1(U_t,\theta)= \frac{\partial \log f(U_t, \theta)}{\partial \theta_k }. \frac{\partial \log f(U_t, \theta)}{\partial \theta_k }+ \frac{\partial^2 \log f(U_t, \theta)}{\partial \theta_k^2}= \frac{\left(\frac{\textbf{I}_{x_i=k}}{n_i}\right)^2}{\theta_k^2}+\frac{-\left(\frac{\textbf{I}_{x_i=k}}{n_i}\right)}{\theta_k^2}$$

$$d_2(U_t,\theta)=\frac{\partial \log f(U_t, \theta)}{\partial \theta_k }. \frac{\partial \log f(U_t, \theta)}{\partial \theta_l }+ \frac{\partial^2 \log f(U_t, \theta)}{\partial \theta_k\partial \theta_l}=\frac{\left(\frac{\textbf{I}_{x_i=k}}{n_i}\right)\left(\frac{\textbf{I}_{x_j=l}}{n_j}\right)}{\theta_k\theta_l} $$

$$d_3(U_t,\theta)=\frac{\partial \log f(U_t, \theta)}{\partial \theta_k }. \frac{\partial \log f(U_t, \theta)}{\partial \eta_{kl} }+ \frac{\partial^2 \log f(U_t, \theta)}{\partial \theta_k\partial \eta_{kl}}=\frac{\left(\frac{\textbf{I}_{x_i=k}}{n_i}\right)}{\theta_k}\left( \frac{y_{ij}}{\eta_{kl}}-\frac{\left(1 -y_{ij} \right)}{1-\eta_{kl}}\right) $$

$$d_4(U_t,\theta)=\frac{\partial \log f(U_t, \theta)}{\partial \theta_l }. \frac{\partial \log f(U_t, \theta)}{\partial \theta_l }+ \frac{\partial^2 \log f(U_t, \theta)}{\partial \theta_l^2}=\frac{\left(\frac{\textbf{I}_{x_j=l}}{n_j}\right)^2}{\theta_l^2} +  \frac{-\left(\frac{\textbf{I}_{x_j=l}}{n_j}\right)}{\theta_l^2} $$

$$d_5(U_t,\theta)=\frac{\partial \log f(U_t, \theta)}{\partial \theta_l }. \frac{\partial \log f(U_t, \theta)}{\partial \eta_{kl} }+ \frac{\partial^2 \log f(U_t, \theta)}{\partial \theta_l\partial \eta_{kl}}=\frac{\left(\frac{\textbf{I}_{x_j=l}}{n_j}\right)}{\theta_l}\left( \frac{y_{ij}}{\eta_{kl}}-\frac{\left( 1 -y_{ij} \right)}{1-\eta_{kl}}\right) $$

$$d_6(U_t,\theta)=\frac{\partial \log f(U_t, \theta)}{\partial \eta_{kl} }. \frac{\partial \log f(U_t, \theta)}{\partial \eta_{kl} }+ \frac{\partial^2 \log f(U_t, \theta)}{\partial \eta_{kl}^2}= \left(\frac{y_{ij}}{\eta_{kl}}-\frac{\left( 1 -y_{ij} \right)}{1-\eta_{kl}} \right)^2+ \frac{-y_{ij}}{\eta_{kl}^2}-\frac{ \left( 1-y_{ij} \right)}{(1-\eta_{kl})^2}$$

The test will be based on $\displaystyle D_{n}(\hat{\theta})=\binom{n}{2}^{-1}\sum_{t=1}^{\binom{n}{2 }}d_l(U_t, \hat{\theta})$, which are the elements of $A_n(\hat{\theta})+B_n(\hat{\theta})$.

Let $l=1, \ldots\, q=p(p+1)/2=6$, define the vector $d(U_t,\theta)$, of dimension $q \times 1 = 6 \times 1$ , so
$$d(U_t,\theta)=\left( \begin{matrix}

  \frac{\left(\frac{\textbf{I}_{x_i=k}}{n_i}\right)^2}{\theta_k^2}+\frac{-\left(\frac{\textbf{I}_{x_i=k}}{n_i}\right)}{\theta_k^2} \\ 

\frac{\left(\frac{\textbf{I}_{x_i=k}}{n_i}\right)\left(\frac{\textbf{I}_{x_j=l}}{n_j}\right)}{\theta_k\theta_l} \\

\frac{\left(\frac{\textbf{I}_{x_i=k}}{n_i}\right)}{\theta_k}\left( \frac{y_{ij}}{\eta_{kl}}-\frac{\left(1 -y_{ij} \right)}{1-\eta_{kl}}\right) \\

\frac{\left(\frac{\textbf{I}_{x_j=l}}{n_j}\right)^2}{\theta_l^2} +  \frac{-\left(\frac{\textbf{I}_{x_j=l}}{n_j}\right)}{\theta_l^2} \\

\frac{\left(\frac{\textbf{I}_{x_j=l}}{n_j}\right)}{\theta_l}\left( \frac{y_{ij}}{\eta_{kl}}-\frac{\left( 1 -y_{ij} \right)}{1-\eta_{kl}}\right) \\

 \left(\frac{y_{ij}}{\eta_{kl}}-\frac{\left( 1 -y_{ij} \right)}{1-\eta_{kl}} \right)^2+ \frac{-y_{ij}}{\eta_{kl}^2}-\frac{ \left( 1-y_{ij} \right)}{(1-\eta_{kl})^2}

\end{matrix}\right)$$

So $D_n(\hat{\theta})=\binom{n}{2}^{-1}\sum_{t=1}^{\binom{n}{2}}d(U_t, \hat{\theta})$. For SBM,

$$D_n(\hat{\theta_{n}})=\left( \begin{matrix}

 \binom{n}{2}^{-1}\sum_{t=1}^{\binom{n}{2}}\left( \frac{\left(\frac{\textbf{I}_{x_i=k}}{n_i}\right)^2}{\hat{\theta_k}^2}+\frac{-\left(\frac{\textbf{I}_{x_i=k}}{n_i}\right)}{\hat{\theta_k}^2} \right) \\ 

\binom{n}{2}^{-1}\sum_{t=1}^{\binom{n}{2}}\left(\frac{\left(\frac{\textbf{I}_{x_i=k}}{n_i}\right)\left(\frac{\textbf{I}_{x_j=l}}{n_j}\right)}{\hat{\theta_k}\hat{\theta_l}} \right)\\

\binom{n}{2}^{-1}\sum_{t=1}^{\binom{n}{2}}\left(\frac{\left(\frac{\textbf{I}_{x_i=k}}{n_i}\right)}{\hat{\theta_k}}\left( \frac{y_{ij}}{\hat{\eta_{kl}}}-\frac{\left(1 -y_{ij} \right)}{1-\hat{\eta_{kl}}}\right)\right) \\

\binom{n}{2}^{-1}\sum_{t=1}^{\binom{n}{2}}\left(\frac{\left(\frac{\textbf{I}_{x_j=l}}{n_j}\right)^2}{\hat{\theta_l}^2} +  \frac{-\left(\frac{\textbf{I}_{x_j=l}}{n_j}\right)}{\hat{\theta_l}^2} \right)\\

\binom{n}{2}^{-1}\sum_{t=1}^{\binom{n}{2}}\left(\frac{\left(\frac{\textbf{I}_{x_j=l}}{n_j}\right)}{\hat{\theta_l}}\left( \frac{y_{ij}}{\hat{\eta_{kl}}}-\frac{\left( 1 -y_{ij} \right)}{1-\hat{\eta_{kl}}}\right)\right) \\

\binom{n}{2}^{-1}\sum_{t=1}^{\binom{n}{2}}\left( \left(\frac{y_{ij}}{\hat{\eta_{kl}}}-\frac{\left( 1 -y_{ij} \right)}{1-\hat{\eta_{kl}}} \right)^2+ \frac{-y_{ij}}{\hat{\eta_{kl}}^2}-\frac{ \left( 1-y_{ij} \right)}{(1-\hat{\eta_{kl}})^2}\right)

\end{matrix}\right)$$

 From $D_n(\hat{\theta})$, we define $\nabla D_n(\hat{\theta})=\binom{n}{2}^{-1}\sum_{t=1}^{\binom{n}{2}}\left[\frac{\partial d_l(U_t, \hat{\theta}) }{ \partial \theta_k}\right].$

$$\nabla D_n(\hat{\theta})= \left( \begin{matrix}
\binom{n}{2}^{-1}\sum_{t=1}^{\binom{n}{2}}\left(\frac{\partial d_1(u, \hat{\theta}) }{ \partial \theta_k}\right) & \binom{n}{2}^{-1}\sum_{t=1}^{\binom{n}{2}}\left(\frac{\partial d_1(u, \hat{\theta}) }{ \partial \theta_l}\right)= 0 & \binom{n}{2}^{-1}\sum_{t=1}^{\binom{n}{2}}\left(\frac{\partial d_1(u, \hat{\theta}) }{ \partial \eta_{kl}}\right)=0 \\ 

\binom{n}{2}^{-1}\sum_{t=1}^{\binom{n}{2}}\left(\frac{\partial d_2(u, \hat{\theta}) }{ \partial \theta_k}\right) & \binom{n}{2}^{-1}\sum_{t=1}^{\binom{n}{2}}\left(\frac{\partial d_2(u, \hat{\theta}) }{ \partial \theta_l}\right) & \binom{n}{2}^{-1}\sum_{t=1}^{\binom{n}{2}}\left(\frac{\partial d_2(u, \hat{\theta}) }{ \partial \eta_{kl}}\right)=0 \\  

\binom{n}{2}^{-1}\sum_{t=1}^{\binom{n}{2}}\left(\frac{\partial d_3(u, \hat{\theta}) }{ \partial \theta_k}\right) & \binom{n}{2}^{-1}\sum_{t=1}^{\binom{n}{2}}\left(\frac{\partial d_3(u, \hat{\theta}) }{ \partial \theta_l}\right)= 0 & \binom{n}{2}^{-1}\sum_{t=1}^{\binom{n}{2}}\left(\frac{\partial d_3(u, \hat{\theta}) }{ \partial \eta_{kl}}\right)\\

\binom{n}{2}^{-1}\sum_{t=1}^{\binom{n}{2}}\left(\frac{\partial d_4(u, \hat{\theta}) }{ \partial \theta_k}\right)=0 & \binom{n}{2}^{-1}\sum_{t=1}^{\binom{n}{2}}\left(\frac{\partial d_4(u, \hat{\theta}) }{ \partial \theta_l}\right) & \binom{n}{2}^{-1}\sum_{t=1}^{\binom{n}{2}}\left(\frac{\partial d_4(u, \hat{\theta}) }{ \partial \eta_{kl}}\right)=0 \\ 

\binom{n}{2}^{-1}\sum_{t=1}^{\binom{n}{2}}\left(\frac{\partial d_5(u, \hat{\theta}) }{ \partial \theta_k}\right)=0 & \binom{n}{2}^{-1}\sum_{t=1}^{\binom{n}{2}}\left(\frac{\partial d_5(u, \hat{\theta}) }{ \partial \theta_l}\right) & \binom{n}{2}^{-1}\sum_{t=1}^{\binom{n}{2}}\left(\frac{\partial d_5(u, \hat{\theta}) }{ \partial \eta_{kl}}\right) \\

\binom{n}{2}^{-1}\sum_{t=1}^{\binom{n}{2}}\left(\frac{\partial d_6(u, \hat{\theta}) }{ \partial \theta_k}\right)=0 & \binom{n}{2}^{-1}\sum_{t=1}^{\binom{n}{2}}\left(\frac{\partial d_6(u, \hat{\theta}) }{ \partial \theta_l}\right)= 0 & \binom{n}{2}^{-1}\sum_{t=1}^{\binom{n}{2}}\left(\frac{\partial d_6(u, \hat{\theta}) }{ \partial \eta_{kl}}\right) 

 \end{matrix} \right)
 $$

The functions that must be summed in the components of $\nabla D_n(\theta)$ are the functions that were obtained in the assumption \ref{supo8}.

$V(\theta_{*})$ is the asymptotic covariance matrix of $\sqrt{n}D_n(\hat{\theta})$ and we have a consistent estimator for $V(\theta_{*})$ It's
{\small $$V_n(\hat{\theta})=\binom{n}{2}^{-1}\sum_{t=1}^{\binom{n}{2}}\left[d (U_t,\hat{\theta})-\nabla D_n(\hat{\theta})A_n(\hat{\theta})^{-1}\nabla \log f(U_t,\hat{\theta} )\right]. \\ \left[d(U_t,\hat{\theta})-\nabla D_n(\hat{\theta})A_n(\hat{\theta})^{-1}\nabla \log f(U_t, \hat{\theta})\right]'.$$}

In SBM, for each of the $U_t$ vectors we will get
$$V_{t}(\hat{\theta})=\left[d(U_t,\hat{\theta})-\nabla D_n(\hat{\theta})A_n(\hat{\theta}) ^{-1}\nabla \log f(U_t,\hat{\theta})\right]. \\ \left[d(U_t,\hat{\theta})-\nabla D_n(\hat{\theta})A_n(\hat{\theta})^{-1}\nabla \log f(U_t, \hat{\theta})\right]',$$
$V_{t}(\hat{\theta})$ is a matrix of dimension $6 \times 6$, and $V_n(\hat{\theta})$ will be a matrix $6 \times 6$, whose entries are the average of the respective entries of the $V_{t}(\hat{\theta})$.

Thus, we arrive at the following test:

\begin{teo}\textbf{(Misspecification Test on Stochachistic Block Models)} If the assumptions from \ref{supo1} to \ref{supo10} are satisfied, and if $h(U)=f(U, \theta_0)$, to $\theta_0 \in \Theta$, then
\begin{equation}
\label{tmge}\mathcal{I}_n=nD_n(\hat{\theta})'(V_n(\hat{\theta}))^{-1}D_n(\hat{\theta})
\end{equation}
has asymptotic distribution $\chi^2_6$.
\label{teo7}
\end{teo}

For a hypothesis test with $\alpha$ of significance, calculate \ref{teo7} and use the criterion: if $\mathcal{I}_n \leq \chi^2_{(\alpha, 6) } $, where $\chi^2_{(\alpha, 6)}$ is the cumulative distribution value of $\chi^2_{6}$ in $\alpha$, the model was well specified, otherwise the model was misspecified.

\section{Simulation}

In the previous sections we built misspecification tests for the ERG and the SBM, specifying the test statistics and their respective distributions. In this Section we will present simulations to verify the behavior of the proposed tests in known scenarios.

The tests were applied to samples of graphs that were generated by varying the values of each parameter of the models, the number of vertices and number of classes, in the case of SBM. To apply the tests, we simulated random samples of graphs in two scenarios: $(1)$ graphs generated from each of the models and $(2)$ graphs generated with other models. From the samples, we estimated the parameters via maximum likelihood and built the sample vectors $U_t$, which will be used in each test. Then we find the auxiliary matrices $d(U_t, \hat{\theta})$, $D_n(\hat{\theta})$, $\nabla D_n(\hat{\theta})$, $A_n( \hat{\theta})$, $\nabla \log f(U_t, \hat{\theta})$ to calculate $V_n(\hat{\theta})$ and finally get $\mathcal{I} _n$ and compare with the cumulative distribution value of $\chi^2_{q}$, fixed to a $\alpha$.

The codes were implemented in R, version $4.1.2$, and executed on a computer
with AMD Ryzen processor $5$ $3.6$ GHz, $16$GB of RAM and $512GB$ SSDM2. The simulation time in each scenario was directly related to the number of vertices and number of classes, in the case of SBM. Some tests ran in seconds while others took a few minutes.

For each of the scenarios, we performed Monte Carlo(MC) replications and investigated what percentage of the tests indicated when the model was well-specified. In Tables \ref{tab1}, \ref{tab2}, \ref{tab3} and \ref{tab4} it is possible to observe the proportion of simulated tests that indicate a well-specified model in each of the scenarios.

\subsection{Behavior of Misspecification Tests on ERG}

For the ERG, we consider tests in which the graphs are divided into two scenarios. In scenario 1, we generate graphs with ERG distribution with only one parameter $\theta$, in scenario 2, we will generate graphs that, intentionally, do not have ERG distribution with only one parameter $\theta$. In both scenarios we varied the number of vertices of the graphs, did $1000$ or $10000$ Monte Carlo(MC) replications and investigated what percentage of the tests indicated when the model was well-specified.

\begin{enumerate}

\item {\bf Scenario 1}: In each replication we generate a random $\alpha$ in the range $(0,1)$, the graphs have $n=50$, $n=100$, $n=200$, $n=1000$, or $n=10000$ vertices. The probability of a vertex $i$ connecting to any other vertex $j$ is exactly $\alpha$ for any vertices $i$ and $j$, that is, theoretically, we are in a graph with an ERG distribution with only one parameter $\theta= \log\left( \frac{\alpha}{1-\alpha}\right)$.

\item {\bf Scenario 2}: In each replication we generate a random $\alpha$ in the range $(0,1)$, the graphs have $n=50$, $n=100$, $n=200$, $n=1000$ or $n=10000$ vertices. The probability of a vertex $i$ connecting to another vertex $j$ is $p_k . \alpha$ for each $\frac{n}{10}$ vertices, where $p_k$ is a random value in the range $(0,1)$, for $k=1,2, \ldots, 10$, that is, we generate graphs whose probability for each group of $\frac{n}{10}$ vertices is a random percentage of $\alpha$, so, theoretically, we are in a graph whose distribution is not that of the ERG with parameter $ \theta= \log\left( \frac{\alpha}{1-\alpha}\right)$.

\end{enumerate}

The test for misspecification of the ERG model has the following assumptions:

\bigskip

$H_0:$ the generating graph has ERG distribution(${\theta}$),
\\against

$H_a:$ the generating graph has no distribution ERG(${\theta}$).





\bigskip

Let's take as an estimate of $\theta$ the maximum likelihood estimator given by
$$ \hat{\theta}=\log \left( \frac{\frac{\sum_{t=1}^{\binom{n}{2}}U_t}{\binom{n}{2}} }{1-\frac{\sum_{t=1}^{\binom{n}{2}}U_t}{\binom{n}{2}}}\right).$$

In both scenarios, following the process described in section $3.1$, we find the sample vectors $U_t$ and the auxiliary matrices $d(U_t, \hat{\theta})$, $D_n(\hat{\theta})$ , $\nabla D_n(\hat{\theta})$, $A_n(\hat{\theta})$, $\nabla \log f(U_t, \hat{\theta})$ to calculate $V_n(\hat{\theta})$. Finally, we get $\mathcal{I}_n$ and compare it to the cumulative distribution value of $\chi^2_{1}$ at $\alpha=0.05$. If $\mathcal{I}_n$ is less than or equal to the value of the cumulative distribution of $\chi^2_{1}$ in $\alpha=0.05$, we do not reject $H_0$, that is, the test of hypotheses indicates that the model was well-specified.

Table \ref{tab1} presents the proportion of simulated tests, in scenario 1, that indicate a well-specified model, we can infer that for different values of $\theta$ in graphs with $50$, $100$, $200$, $1000$ and $10000$ vertices the hypothesis test with a significance level of $5\%$ indicated good adequacy, as expected.

\begin{table}[H]
\centering
\caption{Scenario 1 for the ERG: Proportion of simulated hypothesis tests that indicate a well-specified model.}
\vspace{0.5cm}
\begin{tabular}{|c|c|c|}
\hline 
\text{ Replications} & \text{Number of vertices} & \text{ Proportion } \\
\hline                             
1000 & 50 & 0.998 \\
\hline
1000 & 100 & 0.999 \\
\hline
1000 & 200 & 0.999 \\
\hline
1000 & 1000 & 1 \\
\hline
1000 & 10000 & 1 \\
\hline
10000 & 50 & 0.9984\\
\hline
10000 & 100 & 0.9993 \\
\hline
10000 & 200 & 0.9999\\
\hline
10000 & 1000 &  1 \\
\hline
10000 & 10000 & 1 \\
\hline
\end{tabular}
\label{tab1}
\end{table}

Table \ref{tab2} shows the proportion of simulated tests, in scenario 2, that indicate a well-specified model, we can infer that for different values of $\theta$ in graphs with $50$, $100$, $200$, $1000$ and $10000$ vertices the hypothesis test with a significance level of $5\%$ indicated poor adequacy, as was also expected.

\begin{table}[H]
\centering
\caption{Scenario 2 for the ERG: Proportion of simulated hypothesis tests that indicate a well-specified model.}
\vspace{0.5cm}
\begin{tabular}{|c|c|c|}
\hline 
\text{ Replications} & \text{Number of vertices} & \text{ Proportion } \\
\hline                             
1000 & 50 & 0 \\
\hline
1000 & 100 & 0 \\
\hline
1000 & 200 & 0 \\
\hline
1000 & 1000 & 0\\
\hline
1000 & 10000 &  0\\
\hline
10000 & 50 & 0.0029\\
\hline
10000 & 100 & 0.0014 \\
\hline
10000 & 200 & 0.0006 \\
\hline
10000 & 1000 &  0\\
\hline
10000 & 10000 & 0\\
\hline
\end{tabular}
\label{tab2}
\end{table}

\newpage

\subsection{Behavior of Misspecification Tests on SBM}

For SBM we also consider tests in which the graphs are divided into two scenarios: in {\bf Scenario 1}, we generate the graphs with the SBM distribution with the parameters $\left(\theta, \eta \right)$. In {\bf Scenario 2}, we generate random graphs by perturbing the original parameters $\left(\theta, \eta \right)$. In both scenarios, we varied the number of vertices and the number of blocks of the graphs, did $100$ or $1000$ Monte Carlo(MC) replications and investigated what percentage of the tests indicated when the model was well-specified.
\begin{enumerate}

\item {\bf Scenario 1}: In each replication we generate the matrix $\eta$ where $\eta_{kl}$ is random in $(0,1)$ and $\eta_{kl}=\eta_{lk }$, for $k=1,2, \ldots, m$ and $l=1,2, \ldots, m$, where $m$ is the number of existing blocks and we generate the vector $X$ such that , for $i=1,2, \ldots, n$, $x_i$ is random between $\left\{1,2, \ldots , m \right\}$ and represents the class of vertex $i$. The probability of a vertex $i$ connecting to any other vertex $j$ is exactly $\eta_{kl}$ if vertex $i$ is of class $k$ and vertex $j$ is of class $l$ , that is, a graph with SBM distribution with parameter $\left(\theta, \eta \right)$ where $\eta$ is formed by $\eta_{kl}$ and $\theta_k=\frac{\sum_{i=1}^n{\textit{I}_{x_i=k}}}{m}$.

We generate graphs with $n=90$ and $m=3$ or $m=6$, $n=120$ and $m=4$ or $m=6$, $n=200$ and $m=4$ or $m=10$, $n=300$ and $m=3$ or $m=10$ or $m=15$, $n=1000$ and $m=4$ or $m=10$ or $ m=20$ or $m=100$.

\item {\bf Scenario 2}: In each replication we generate the matrix $\eta$ where $\eta_{kl}$ is random in $(0,1)$ and $\eta_{kl}=\eta_{lk }$, for $k=1,2, \ldots, m$ and $l=1,2, \ldots, m$, where $m$ is the number of existing blocks and we generate the vector $X$ such that , for $i=1,2, \ldots, n$, $x_i$ is random between $\left\{1,2, \ldots , m \right\}$ and represents the class of vertex $i$. The probability of a vertex $i$ connecting to another vertex $j$ is $p_k . \eta$ for each $\frac{n}{10}$ vertices, where $p_k$ is a random value in the range $(0,1)$, for $k=1,2, \ldots, 10$, that is, we generate graphs whose probability for each group of $\frac{n}{10}$ vertices is a random percentage of $\eta$, thus, a graph whose distribution is not that of the SBM with the parameters $\left( \theta, \eta \right)$ where $\eta$ is formed by $\eta_{kl}$ and $\theta_k=\frac{\sum_{i=1}^n{\textit{I}_{x_i=k}}}{m}$. We generate graphs where $n=90$ and $m=3$ or $m=6$, $n=120$ and $m=4$ or $m=6$, $n=200$ and $m=4 $ or $m=10$ and $n=300$ and $m=3$ or $m=10$ or $m=30$.

\end{enumerate}







The SBM misspecification test has the following assumptions:

\bigskip

$H_0:$ the generating graph has distribution SBM$\left(\theta, \eta \right)$,
\\against

$H_a:$ the generating graph has no distribution SBM$\left(\theta, \eta \right)$.





\bigskip

let's take as an estimate of $\left(\theta, \eta \right)$ the maximum likelihood estimator $\left(\hat{\theta}, \hat{\eta} \right)$ obtained via Variational EM by the package \textit{blockmodels} from R software,

In both scenarios, following the process described in section $3.2$, we construct the sample vectors $U_t$. Then we build the auxiliary matrices $d(U_t, \hat{\theta})$, $D_n(\hat{\theta})$, $\nabla D_n(\hat{\theta})$, $A_n( \hat{\theta})$, $\nabla \log f(U_t, \hat{\theta})$ to calculate $V_n(\hat{\theta})$ and finally get $\mathcal{I} _n$ and compare with the cumulative distribution value of $\chi^2_{6}$ at $\alpha=0.05$. If $\mathcal{I}_n$ is less than or equal to the value of the cumulative distribution of $\chi^2_{6}$ in $\alpha=0.05$, $H_0$ is accepted, that is, the test of hypotheses indicates that the model was well-specified as a SBM.

Table \ref{tab3} presents the proportion of simulated tests, in scenario 1, that indicate a well-specified model, we can infer that for different values of $\left(\theta, \eta \right)$ in graphs with variations of the number of blocks and vertices the test indicated good fit as expected. 

\begin{table}[H]
\centering
\caption{Scenario 1 for the SBM: Proportion of simulated hypothesis tests that indicate a well-specified model.}
\vspace{0.5cm}
\begin{tabular}{|c|c|c|c|}
\hline 
\text{ Replication} & \text{Number of vertices} & \text{Number of blocks} & \text{ Proportion } \\
\hline                             
100 & 90 & 3 & 1\\
\hline
100 & 90 & 6 &  0.99\\
\hline
100 & 120 & 4 & 1\\
\hline
100 & 120 & 6 & 0.98\\
\hline
100 & 200 & 4 &  0.99\\
\hline
100 & 200 & 10 & 0.97 \\
\hline
100 & 300 & 3 & 1\\
\hline
100 & 300 & 10 &  1\\
\hline
100 & 300 & 30 &  0.99\\
\hline
1000 & 90 & 3 & 0.963\\
\hline
1000 & 90 & 6 &  0.972\\
\hline
1000 & 120 & 4 & 0.968\\
\hline
1000 & 120 & 6 & 0.981\\
\hline
1000 & 200 & 4 &  0.971\\
\hline
1000 & 200 & 10 & 0.982\\
\hline
1000 & 300 & 3 & 0.983\\
\hline
1000 & 300 & 10 &  0.989\\
\hline
1000 & 300 & 30 & 0.991 \\
\hline
\end{tabular}
\label{tab3}
\end{table}

Table \ref{tab4} presents the proportion of simulated tests, in scenario 2, that indicate a well-specified model, we can infer that for different values of $\left(\theta, \eta \right)$ in graphs with variations of the number of blocks and vertices, the test indicated a bad fit, as was also expected.

\begin{table}[H]
\centering
\caption{Scenario 2 for the SBM: Proportion of simulated hypothesis tests that indicate a well-specified model.}
\vspace{0.5cm}
\begin{tabular}{|c|c|c|c|}
\hline 
\text{ Replication} & \text{Number of vertices} & \text{Number of blocks} & \text{ Proportion } \\
\hline                             
100 & 90 & 3 & 0\\
\hline
100 & 90 & 6 &  0.02\\
\hline
100 & 120 & 4 & 0\\
\hline
100 & 120 & 6 & 0.01\\
\hline
100 & 200 & 4 &  0\\
\hline
100 & 200 & 10 & 0.02\\
\hline
100 & 300 & 3 & 0.01\\
\hline
100 & 300 & 10 &  0.02\\
\hline
100 & 300 & 30 &  0.02\\
\hline
1000 & 90 & 3 & 0.042\\
\hline
1000 & 90 & 6 &  0.054\\
\hline
1000 & 120 & 4 & 0.036\\
\hline
1000 & 120 & 6 & 0.044\\
\hline
1000 & 200 & 4 & 0.031 \\
\hline
1000 & 200 & 10 & 0.038\\
\hline
1000 & 300 & 3 & 0.022\\
\hline
1000 & 300 & 10 & 0.024 \\
\hline
1000 & 300 & 30 & 0.032 \\
\hline
\end{tabular}
\label{tab4}
\end{table}

\newpage

\subsection{Conclusions}
The simulations were performed with random samples of graphs in two scenarios for each of the analyzed models: ERG and SBM: (1) In the first scenario, the graphs were generated from each of the specified models. (2) In the second scenario, the graphs were generated with other models. Parameter estimates were made via maximum likelihood and we calculated test statistics for each sample. Tables \ref{tab1} and \ref{tab3} refer to Scenario 1, where graphs were generated according to the ERG and SBM models, showing that most of the simulated tests indicate that the model is adequate to the data. But when the generated graphs are not generated following the models, the results shown in Tables \ref{tab2} and \ref{tab4}, referring to scenario 2, show that most of the tests applied indicate that the model is not suitable for the data .

So, we observe that when the test is applied to a sample that, in fact, was generated according to the tested models, it gets it right (does not reject $H_0$) in practically all simulations. For both ERG and SBM , Tables \ref{tab1} and \ref{tab3}.

In the case of the test for specification error of the ERG, we noticed from Table \ref{tab2} that, even though the number of vertices is considerably large, the proportion of errors is very small, less than $1\%$ in all replicas.

We can see from Table \ref{tab4} that, even when we have a large number of vertices and blocks, the error proportion of the proposed test is also very small, lower than the expected $5\%$ of the type I error.

With these simulations, we illustrate the effectiveness of the proposed tests to verify the specification error of the ERG and SBM models.

The R codes used to generate the graphs and calculate the tests can be found at:

\textit{ https://drive.google.com/drive/folders/1epDTdqS42853\_Arrg7TzMsz3oNaLQfbG *}

*select and copy the address, paste in the browser and add the \_ symbol manually.

\end{document}